\documentclass[11pt]{article}
\usepackage{geometry}                
\usepackage{amssymb}
\usepackage{amsmath}
\usepackage{amsfonts}
\usepackage{color}
\usepackage{graphicx}
\usepackage{float}
\usepackage{graphicx}
\usepackage{epstopdf}
\usepackage[round]{natbib}

\usepackage{times}
\usepackage[parfill]{parskip}

\newcommand{\R}{{\rm I\!R}}

\newcommand{\argmax}[1]{\mathop{\hbox{argmax}}_{#1}}

\newcommand{\diag}[1]{\hbox{diag}(#1)}

\newtheorem{theorem}{Theorem}

\newenvironment{problem}[1][Problem]{\begin{trivlist}
\item[\hskip \labelsep {\bfseries #1}]}{\end{trivlist}}

\newcommand{\consto}[1]{\sigma_{#1}(A,2)}
\newcommand{\constobar}[1]{\sigma_{#1}(\bar{A},2)}
\newcommand{\constok}[1]{\sigma_{#1}(A_k,2)}
\newcommand{\constobark}[1]{\sigma_{#1}(\bar{A}_k,2)}

\newcommand{\const}[1]{\sigma_{#1}(A,\infty)}
\newcommand{\constbar}[1]{\sigma_{#1}(\bar{A},\infty)}

\newcommand\numberthis{\addtocounter{equation}{1}\tag{\theequation}}

\def\norm#1{\|#1\|}

\title{Convergence Rates for Greedy Kaczmarz Algorithms, and Faster Randomized Kaczmarz Rules Using the Orthogonality Graph}

\author{
Julie Nutini$^1$, Behrooz Sepehry$^1$, Issam Laradji$^1$, \\Mark Schmidt$^1$, Hoyt Koepke$^2$, and Alim Virani$^1$\\\\
$^1$The University of British Columbia, $^2$Dato}
\date{}

\begin{document}
\maketitle

\begin{abstract}

The Kaczmarz method is an iterative  algorithm for solving systems of linear equalities and inequalities, that iteratively projects onto these constraints. Recently, Strohmer and Vershynin [{\em J.~Fourier Anal. Appl.,} 15(2):262-278, 2009] gave a non-asymptotic convergence rate analysis for this algorithm, spurring numerous extensions and generalizations of the Kaczmarz method. Rather than the randomized selection rule analyzed in that work, in this paper we instead discuss greedy and approximate greedy selection rules. We show that in some applications the computational costs of greedy and random selection are comparable, and that in many cases greedy selection rules give faster convergence rates than random selection rules.  Further, we give the first multi-step analysis of Kaczmarz methods for a particular greedy rule, and propose a provably-faster randomized selection rule for matrices with many pairwise-orthogonal rows.

\end{abstract}

\section{Kaczmarz Method}

Solving large linear systems is a fundamental problem in machine learning. Applications range from least-squares problems to Gaussian processes to graph-based semi-supervised learning. All of these applications (and many others) benefit from advances in solving large-scale linear systems. The Kaczmarz method is a particular iterative algorithm suited for solving consistent linear systems of the form $Ax =b$. It was originally proposed by Polish mathematician Stefan \citet{kaczmarz1937} and later re-invented by \citet{gordon1970ART} under the name {\em algebraic reconstruction technique} (ART).  It has been used in numerous applications including image reconstruction and digital signal processing, and belongs to several general categories of methods including {\em row-action}, {\em component-solution}, {\em cyclic projection}, and {\em successive projection} methods~\citep{censor1981}.

At each iteration $k$, the Kaczmarz method uses a \emph{selection rule} to choose some row $i_k$ of $A$ and then projects the current iterate $x^k$ onto the corresponding hyperplane $a_{i_k}^Tx^k = b_{i_k}$. Classically, the two categories of selection rules are \emph{cyclic} and \emph{random}. Cyclic selection repeatedly cycles through the coordinates in sequential order, making it simple to implement and computationally inexpensive. There are various linear convergence rates for cyclic selection~\citep[see][]{deutsch1985,deutsch1997,galantai2005}, but these rates are in terms of cycles through the entire dataset and involve constants that are not easily interpreted. Further, the performance of cyclic selection worsens if we have an undesirable ordering of the rows of $A$.

Randomized selection has recently become the default selection rule in the literature on Kaczmarz-type methods. Empirically, selecting $i_k$ randomly often performs substantially better in practice than cyclic selection~\citep{cenker1992POCS,herman1993}. Although a number of asymptotic convergence rates for randomized selection have been presented~\citep{whitney1967, tanabe1971, censor1983, hanke1990}, the pivotal theoretical result supporting the use of randomized selection for the Kaczmarz method was given by~\citet{strohmer2009}. They  proved a simple non-asymptotic linear convergence rate (in expectation) in terms of the number of iterations, when rows are selected proportional to their squared norms. This work spurred numerous extensions and generalizations of the randomized Kaczmarz method~\citep{needell2010noisy, leventhal2010constraints, zouzais2012extended, lee2013, liu2014accelerated,ma2015}, including similar rates when we replace the equality constraints with inequality constraints.


Rather than cyclic or randomized, in this work we consider {\em greedy} selection rules. 
There are very few results in the literature that explore the use of greedy selection rules for Kaczmarz-type methods. \citet{griebel2012} present the {\em maximum residual rule} for multiplicative Schwarz methods, for which the randomized Kaczmarz iteration is a special case. Their theoretical results show similar convergence rate estimates for both greedy and random methods, suggesting there is no advantage of greedy selection over randomized selection (since greedy selection has additional computational costs). \citet{eldar2011johnsonlindenstrauss}  propose a greedy {\em maximum distance rule}, which they approximate using the Johnson-Lindenstrauss~\citeyearpar{johnson1984lemma} transform to reduce the computation cost. They show that this leads to a faster algorithm in practice, and show that this rule may achieve more progress than random selection on certain iterations.

In the next section, we define several relevant problems of interest  in machine learning that can be solved via Kaczmarz methods. Subsequently, we define the greedy selection rules and discuss cases where they can be computed efficiently.
In Section \ref{sec:rules} we give faster convergence rate analyses for both the maximum residual rule and the maximum distance rule, which clarify the relationship of these rules to random selection and show that greedy methods will typically have better convergence rates than randomized selection. Section \ref{sec:cd} contrasts Kaczmarz methods with coordinate descent methods, Section~\ref{sect:diagonal} considers a simplified setting where we explicitly compute the constants in the convergence rates,
 Section \ref{sec:approx} considers how these convergence rates are changed under \emph{approximations} to the greedy rules, and Section~\ref{sec:inequalities} discusses the case of inequality constraints. We further give a non-trivial \emph{multi-step} analysis of the maximal residual rule (Section \ref{sec:multi}), which is the first multi-step analysis of any Kaczmarz algorithm. By taking the multi-step perspective, we also propose provably-faster randomized selection rules for matrices $A$ with pairwise-orthogonal rows by using the so-called ``orthogonality graph". Section \ref{sec:numerics} presents numerical experiments evaluating greedy Kaczmarz methods.


\section{Problems of Interest}\label{sec:problems}
We first consider systems of linear equations,
\begin{equation}\label{linearsys}
	Ax = b,
\end{equation}
where $A$ is an 
$m \times n$ matrix 
and $b \in \R^m$. We assume the system is \emph{consistent}, meaning a solution $x^*$ exists. We denote the rows of A by $a_1^\top, \dots, a_m^\top$, where each $a_i \in R^n$, and use $b = (b_1,\dots, b_m)^\top$, where each $b_i \in \R$. One of the most important examples of a consistent linear system, and a fundamental model in machine learning, is the least squares problem,
\[
	\min_{x \in \R^n} ~\frac{1}{2} \| A x - b \|^2.
\]
An appealing way to write a least squares problem as a linear system is to solve the $(n+m)$-variable consistent system~\citep[see also][]{zouzais2012extended}
\[	
	\begin{pmatrix}
		A & - I \\
		\bf{ 0 } & A^T
	\end{pmatrix}
	\begin{pmatrix}
		x \\
		y
	\end{pmatrix}
	= 
	\begin{pmatrix}
		b \\
		0
	\end{pmatrix}.
\]
 Other applications in machine learning that involve solving consistent linear systems include: least-squares support vector machines, Gaussian processes, fitting the final layer of a neural network (using squared-error), graph-based semi-supervised learning or other graph-Laplacian problems~\citep{bengio2006}, and finding the optimal configuration in Gaussian Markov random fields~\citep{rue2005}. 


Kaczmarz methods can also be applied to solve consistent systems of linear \emph{inequalities},
\[
Ax \leq b,
\]
or combinations of linear equalities and inequalities. We believe there is a lot potential to use this application of Kaczmarz methods in machine learning. Indeed, a classic example of solving linear inequalities is finding a linear separator for a binary classification problem. The classic perceptron algorithm is a generalization of the Kaczmarz method, but unlike the classic sublinear rates of perceptron methods~\citep{novikoff1963convergence} we can show a linear rate for the Kaczmarz method.

Kaczmarz methods could also be used to solve the $\ell_1$-regularized robust regression problem,
\[
\min_x f(x) := \norm{Ax - b}_1 + \lambda\norm{x}_1,
\]
for $\lambda \geq 0$. We can formulate finding an $x$ with $f(x) \leq \tau$ for some constant $\tau$ as a set of linear inequalities. By doing a binary search for $\tau$ and using \emph{warm-starting}, this can be substantially faster than existing approaches like stochastic subgradient methods (which have a sublinear convergence rate) or formulating as a linear program (which is not scaleable due to the super-linear cost). The above logic applies to many piecewise-linear problems in machine learning like variants of support vector machines/regression with the $\ell_1$-norm, regression under the $\ell_\infty$-norm, and linear programming relaxations for decoding in graphical models.

\section{Kaczmarz Algorithm and Greedy Selection Rules}

The Kaczmarz algorithm for solving linear systems begins from an initial guess $x^0$, and each iteration $k$ chooses a row $i_k$ and projects the current iterate $x^k$ onto the hyperplane defined by $a_{i_k}^T x^k = b_{i_k}$. This gives the iteration
\begin{equation}\label{iterupdate}
	x^{k+1} = x^k + \frac{b_{i_k} - a_{i_k}^T x^k}{\| a_{i_k} \|^2} a_{i_k},
\end{equation}  
and the algorithm converges to a solution $x^*$ under weak conditions (e.g., each $i$ is visited infinitely often).

We consider two greedy selection rules: the {\em maximum residual} rule and the {\em maximum distance} rule.
The maximum residual (MR) rule selects $i_k$ according to
\begin{equation}\label{maxResRule}
	i_k = \argmax{i} ~|a_i^T x^k - b_i|,
\end{equation}
which is the equation $i_k$ that is `furthest' from being satisfied. The maximum distance (MD) rule selects $i_k$ according to
\begin{equation}\label{maxDistRule}
	 i_k = \argmax{i} ~ \left | \frac{a_i^T x^k - b_i}{\|a_i\|} \right |,
\end{equation}
which is the rule that maximizes the distance between iterations, $\| x^{k+1} - x^k \|$.

\subsection{Efficient Calculations for Sparse $A$}\label{sec:calc}

In general, computing these greedy selection rules exactly is too computationally expensive, but in some applications we can compute them efficiently.
 For example, 
consider a \emph{sparse} $A$ with at most $c$ non-zeros per column and at most $r$ non-zeros per row.
In this setting, we show in Appendix~\ref{app:calc} that both rules can be computed exactly in $O(cr \log m)$ time,
using that projecting onto row $i$ does not change the residual of row $j$ if $a_i$ and $a_j$ do not share a non-zero index. 

The above sparsity condition guarantees that row $i$ is orthogonal to row $j$,
and indeed projecting onto row $i$  will not change the residual of row $j$ under the more general condition that $a_i$ and $a_j$ are orthogonal. Consider what we call the \emph{orthogonality graph}: an undirected graph on $m$ nodes where we place on edge between nodes $i$ and $j$ if $a_i$ is not orthogonal to $a_j$. Given this graph, to update all residuals after we update a row $i$ we only need to update the neighbours of node $i$ in this graph. Even if $A$ is dense ($r = n$ and $c = m$), if the maximum number of neighbours is $g$, then tracking the maximum residual costs $O(gr + g\log(m))$. If $g$ is small, this could still be comparable to the $O(r + \log(m))$ cost of using existing randomized selection strategies.

\subsection{Approximate Calculation}

Many applications, particularly those arising from graphical models with a simple structure, will allow efficient calculation of the greedy rules using the method of the previous section. However, in other applications it will be too inefficient to calculate the greedy rules. Nevertheless,~\citet{eldar2011johnsonlindenstrauss} show that it's possible to efficiently select an $i_k$ that \emph{approximates} the greedy rules by making use of the dimensionality reduction technique of~\citet{johnson1984lemma}. Their experiments show that approximate greedy rules can be sufficiently accurate and that they still outperform random selection. After first analyzing exact greedy rules in the next section, we analyze the effect of using approximate rules in Section \ref{sec:approx}.

\section{Analyzing Selection Rules}\label{sec:rules}

All the convergence rates we discuss use the following relationship between $\| x^{k+1} - x^* \|$ and $\| x^k - x^* \|$:
\begin{align*}
	\|x^{k+1} - x^* \|^2 
	&= \| x^k - x^* \|^2 \!-\! \|x^{k+1} - x^k \|^2 + 2 \underbrace{\langle x^{k+1} - x^*, x^{k+1} \!-\! x^k \rangle}_{(= 0, \text{ by orthogonality)}}.
\end{align*}
Using the definition of $x^{k+1}$ from~\eqref{iterupdate} and simplifying, we obtain for the selected $i_k$ that
\begin{equation}\label{derivationstart}
	\|x^{k+1} - x^* \|^2 = \| x^k - x^* \|^2 - \frac{\left (a_{i_k}^Tx^k - b_{i_k} \right)^2}{\| a_{i_k} \|^2}.
\end{equation}

\subsection{Randomized and Maximum Residual}\label{subsec:UandMR}

We first give an analysis of the Kaczmarz method with \emph{uniform} random selection of the row to update $i$ (which we abbreviate as `U'). Conditioning on the $\sigma$-field $\mathcal{F}_{k}$ generated by the sequence $\{ x^0, x^1, \dots, x^k \}$, and taking expectations of both sides of \eqref{derivationstart}, when $i_k$ is selected using U we obtain
\begin{align*}
	\mathbb{E} [\| x^{k+1} - x^* \|^2] 
	&~= \| x^k - x^* \|^2 - \mathbb{E} \left [ \frac{\left (a_i^Tx^k - b_i \right)^2}{\| a_i \|^2} \right] \\
	&~= \| x^k - x^* \|^2 - \sum_{i = 1}^m \frac{1}{m} \frac{(a_i^\top (x^k - x^*))^2}{\| a_i \|^2} \\
	&~\le \| x^k - x^* \|^2 - \frac{1}{m \| A \|^2_{\infty,2}}\sum_{i = 1}^m (a_i^\top (x^k - x^*))^2 \\
	&~= \| x^k - x^* \|^2 - \frac{1}{m \| A \|^2_{\infty,2}}\norm{A(x^k - x^*)}^2 \\
	&~\le \left ( 1 - \frac{\consto{}^2}{m \| A \|^2_{\infty,2}} \right ) \| x^k - x^* \|^2 \numberthis \label{uniform},
\end{align*}
where $\| A \|_{\infty,2}^2 := \max_i \{ \| a_i \|^2 \}$ and $\sigma(A,2)$ is the Hoffman~\citeyearpar{hoffman1952} constant. We've assumed that $x^k$ is not a solution, allowing us to use Hoffman's bound. When $A$ has independent columns,
$\sigma(A,2)$ is the $n$th singular value of $A$ and in general it is the smallest non-zero singular value.

The argument above is related to the analysis of~\citet{vishnoi2012laplacian} but is simpler due to the use of the Hoffman bound. Further, this simple argument
makes it straightforward to derive bounds on other rules.
For example, we can derive the convergence rate bound of~\citet{strohmer2009} by following the above steps
 but selecting $i$ non-uniformly with probability $\| a_i \|^2/\| A \|_F^2$ (where $\| A \|_F$ is the Frobenius norm of $A$). We review these steps in Appendix~\ref{app:UandMR}, showing that this non-uniform (NU) selection strategy has
\begin{equation}\label{nonuniform}
	\mathbb{E}[\| x^{k+1} - x^* \|^2] \le \bigg ( 1 - \frac{\consto{}^2}{\| A \|^2_F} \bigg ) \| x^k - x^* \|^2.
\end{equation}
 This strategy requires prior knowledge of the row norms of $A$, but this is a one-time computation and can be reused for any linear system involving $A$.
Because $\norm{A}_F^2 \leq m\norm{A}_{\infty,2}^2$, the NU rate~\eqref{nonuniform} is at least as fast as the uniform rate~\eqref{uniform}.

While a trivial analysis shows that the MR rule also satisfies~\eqref{uniform} in a deterministic sense, in Appendix~\ref{app:UandMR}
we give a tighter analysis of the MR rule showing it has the convergence rate
\begin{equation}\label{maxres}
	\| x^{k+1} - x^* \|^2 \le \bigg ( 1 - \frac{\const{}^2}{\|A\|^2_{\infty,2}} \bigg ) \| x^k - x^* \|^2,
\end{equation}
where the Hoffman-like constant $\const{}$ satisfies the relationship
\[
	\frac{\consto{}}{\sqrt{m}} \le \const{} \le \consto{}.
\]
Thus, at one extreme the maximum residual rule obtains the same rate as \eqref{uniform} for uniform selection when $\consto{}^2/m \approx \const{}^2$. However, at the other extreme the maximum residual rule could be faster than uniform selection by a factor of $m$ ($\const{}^2 \approx \consto{}^2 $). Thus, although the uniform and MR bounds are the same in the worst case, the MR rule can be superior by a large margin.

In contrast to comparing U and MR, the MR rate may be faster or slower than the NU rate. This is because
\[
	\|A \|_{\infty,2} \le \|A \|_{F} \le \sqrt{m} \|A \|_{\infty,2},
\]
so these quantities and the relationship between $\consto{}$ and $\const{}$ influence which bound is tighter.

\subsection{Tighter Uniform and MR Analysis}\label{subsect:tight}

In our derivations of rates \eqref{uniform} and \eqref{maxres}, we use the  inequality
\begin{equation}\label{extraineq}
	\| a_i \|^2 \le  \| A \|^2_{\infty,2} \quad \forall ~i,
\end{equation}
which leads to a simple result but could be very loose if the range of the row norms is large.  In this section, we give tighter analyses of the U and MR rules that are less interpretable but are tighter because they avoid this inequality.

In order to avoid using this inequality for our analysis of U, we can absorb the row norms of $A$ into a row weighting matrix $D$, where $D = \diag{ \|a_1\|, \|a_2\|, \dots, \|a_m \|}$. Defining $\bar{A} := D^{-1}A$, we show in Appendix~\ref{app:tight} that this results in the following upper bound on the convergence rate for uniform random selection,
\begin{equation}\label{uniformtight}
	\mathbb{E} [\| x^{k+1} - x^* \|^2] \le \left ( 1 - \frac{\constobar{}^2}{m} \right ) \| x^k - x^* \|^2.
\end{equation}
A similar result is given by \citet{needell2014sgd} under the stronger assumption that $A$ has independent columns.
The rate in \eqref{uniformtight} is tighter than \eqref{uniform}, since $\consto{}/\| A \|_{\infty,2} \le  \constobar{}$~\citep{vandersluis1969}. Further, \emph{this rate  can be faster than the non-uniform sampling method} of \citet{strohmer2009}. 
For example, suppose row $i$ is orthogonal to all other rows but has a significantly larger row norm than all other row norms. In other words, $\| a_i \| >> \| a_j \|$ for all $j \not = i$. In this case, NU selection will repeatedly select row $i$ (even though it only needs to be selected once), whereas U will only select it on each iteration with probability $1/m$. It has been previously pointed out that Strohmer and Vershynin's method can perform poorly if you have a problem where one row norm is significantly larger than the other row norms~\citep{censor2009notesSV}. This result theoretically shows that U can have a tighter bound than the NU method of Strohmer and Vershynin.

In Appendix~\ref{app:tight}, we also give a simple modification of our analysis of the MR rule, which leads to the rate
\begin{equation}\label{maxrestight}
	\| x^{k+1} - x^* \|^2 \le  \bigg ( 1 - \frac{\const{}^2}{\|a_{i_k}\|^2} \bigg ) \| x^k - x^* \|^2.
\end{equation}
 This bound depends on the \emph{specific} $\|a_{i_k} \|$ corresponding to the $i_k$ selected at each iteration $k$. This convergence rate will be faster whenever we select an ${i_k}$ with $\|a_{i_k} \| < \| A \|_{\infty,2}$. However, in the worst case we repeatedly select $i_k$ values with $\|a_{i_k} \| = \| A \|_{\infty,2}$ so there is no improvement. In Section \ref{sec:multi}, we return to this issue and give tighter bounds on the \emph{sequence} of $\| a_{i_k} \|$ values for problems with sparse orthogonality graphs.

\subsection{Maximum Distance Rule}\label{subsect:maxDist}

If we can only perform one iteration of the Kaczmarz method, the \emph{optimal} rule (with respect to iteration progress) is in fact the MD rule.
In Appendix~\ref{app:maxDist}, we show that this strategy achieves a rate of
\begin{equation}\label{maxdist}
	\| x^{k+1} - x^* \|^2 \le \bigg ( 1 - \constbar{}^2 \bigg ) \| x^k - x^* \|^2,
\end{equation}
where $\constbar{}$ satisfies
\[
\hspace{-6pt}	\max \!\left \{\!\frac{\constobar{}}{\sqrt{m}}, \frac{\consto{}}{\|A\|_F}, \frac{\const{}}{\|A\|_{\infty,2}} \!\right \} \!\le \! \constbar{} \! \le \! \constobar{}.
\]
Thus, the maximum distance rule is at least as fast as the fastest among the U/NU/MR$_\infty$ rules, where MR$_\infty$ refers to rate \eqref{maxres}. Further, in Appendix~\ref{app:73} we show that  this new rate is not only simpler but is strictly tighter than the rate reported by~\citet{eldar2011johnsonlindenstrauss} for the exact MD rule.
%
\begin{table}[h!]
\centering
\caption{Comparison of Convergence Rates} \bigskip
\begin{tabular}{| l | c | c | c | c | c | c |}
\hline
  &U$_\infty$ & U & NU & MR$_\infty$ & MR & MD \\ \hline
  U$_\infty$ & $=$ & $\le$ & $\le$ & $\le$ & $\le$ & $\le$ \\ \hline
  U &  & $=$ & P & P & P & $\le$ \\ \hline
  NU &  &  & $=$ & P & P & $\le$ \\ \hline
  MR$_\infty$ &  &  &  & $=$ & $\le$ & $\le$ \\ \hline
  MR &  &  &  &  & $=$ & $\le$\\ \hline
  MD &  &  &  &  &  & $=$\\
\hline
\end{tabular}
\label{table:1}
\end{table}
In Table~\ref{table:1}, we summarize the relationships we have discussed in this section among the different selection rules.
We use the following abbreviations: U$_\infty$ - uniform \eqref{uniform}, U - tight uniform \eqref{uniformtight}, NU - non-uniform \eqref{nonuniform}, MR$_\infty$ - maximum residual \eqref{maxres}, MR - tight maximum residual \eqref{maxrestight} and MD - maximum distance \eqref{maxdist}. The inequality sign ($\leq$) indicates that the bound for the selection rule listed in the row is slower or equal to the rule listed in the column, while we have written `P' to indicate that the faster method is problem-dependent. 


\section{Kaczmarz and Coordinate Descent}\label{sec:cd}

With the exception of the tighter U and MR rate, the results of the previous section are analogous to the recent results of \citet{nutini2015} for coordinate descent methods. Indeed, if we apply coordinate descent methods to minimize the squared error between $Ax$ and $b$ then we obtain similar-looking rates and analogous conclusions. With cyclic selection this is called the Gauss-Seidel method, and as discussed by~\citet{ma2015} there are several connections/differences between this method and Kaczmarz methods.
In this section we highlight some key differences.

First, the previous work required strong-convexity which would require that $A$ has independent columns. This is often unrealistic, and our results from the previous section hold for any $A$.\footnote{\citet{karimi2016} recently showed that the results of~\citet{nutini2015} apply for general least squares problems.} Second, here our results are in terms of the iterates $\|x^k - x^*\|$, which is the natural measure for linear systems. The coordinate descent results are in terms of $f(x^k) - f(x^*)$ and although it's possible to use strong-convexity to turn this into a rate on $\|x^k - x^*\|$, this would result in a looser bound and would again require strong-convexity to hold~\citep[see][]{ma2015}. On the other hand, coordinate descent gives the least squares solution for inconsistent systems. However, this is also true of the Kaczmarz method using the formulation in Section~\ref{sec:problems}. Another subtle issue is that the Kaczmarz rates depend on the row norms of $A$ while the coordinate descent rates depend on the column norms. Thus, there are scenarios where we expect Kaczmarz methods to be much faster and vice versa. Finally, we note that Kaczmarz methods can be extended to allow inequality constraints (see Section~\ref{sec:inequalities}).

As discussed by~\citet{wright2015}, Kaczmarz methods can also be interpreted as coordinate descent methods on the dual problem
\begin{equation}\label{dual}
	\min_y \frac{1}{2} \| A^Ty \|^2 - b^Ty,
\end{equation}
where $x = A^Ty^*$ so that $Ax = A A^T y^* = b$. Applying the Gauss-Southwell rule in this setting yields the MR rule while applying the Gauss-Southwell-Lipschitz rule yields the MD rule (see Appendix~\ref{app:cd} for details and numerical comparisons, indicating that in some cases Kaczmarz substantially outperforms CD). However, applying the analysis of~\citet{nutini2015} to this dual problem would require that $A$ has independent rows and would only yield a rate on the dual objective, unlike the convergence rates in terms of $\norm{x^k - x^*}$ that hold for general $A$ from the previous section.

\section{Example: Diagonal $A$}\label{sect:diagonal}

To give a concrete example of these rates, we consider the simple case of a diagonal $A$. While such problems are not particularly interesting, this case provides a simple setting to understand these different rates without referring to Hoffman bounds.

Consider a square diagonal matrix $A$ with $a_{ii} > 0$ for all $i$. In this case, the diagonal entries are the eigenvalues $\lambda_i$ of the linear system. The convergence rate constants for this scenario are given in Table~\ref{table:2}.
\bgroup
\def\arraystretch{2}
\begin{table}[ht!]
\caption{Convergence Rate Constants for Diagonal $A$} \bigskip
\centering
\begin{tabular}{| c | c |}
\hline
U$_\infty$ & $\displaystyle \left ( 1 -  \frac{\lambda_m^2}{m\lambda_{1}^2} \right )$ \\ \hline
U & $\left ( 1 -  \displaystyle \frac{1}{m} \right )$ \\ \hline
NU & $\left ( 1 -  \displaystyle \frac{\lambda_m^2}{\sum_i \lambda_i^2} \right )$ \\ \hline
MR$_\infty$& $\left ( 1 -  \displaystyle  \frac{1}{\lambda_1^2} \left [ \sum_i\frac{1}{\lambda_i^2} \right ]^{-1} \right )$ \\ \hline
MR& $\left ( 1 -  \displaystyle  \frac{1}{\lambda_{i_k}^2} \left [ \sum_i\frac{1}{\lambda_i^2} \right ]^{-1} \right )$ \\ \hline
MD & $\left ( 1 -  \displaystyle \frac{1}{m} \right )$ \\
\hline
\end{tabular}
\label{table:2}
\end{table}
We provide the details in Appendix~\ref{app:diagonal} of the derivations for $\const{}$ and $\constbar{}$, as well as substitutions for the uniform, non-uniform, and uniform tight rates to yield the above table. We note that the uniform tight rate follows from $\lambda_{m}^2(\bar{A})$ being equivalent to the minimum eigenvalue of the identity matrix.

If we consider the most basic case when all the eigenvalues of $A$ are equal, then all the selection rules yield the same rate of $(1 - 1/m)$ and the method converges in at most $m$ steps for  greedy selection rules and in at most $O(m \log m)$ steps (in expectation) for the random rules (due to the `coupon collector' problem). Further, this is the worst situation for the greedy MR and MD rules since they satisfy their lower bounds on $\const{}$ and $\constbar{}$.

Now consider the extreme case when all the eigenvalues are equal except for one. For example, consider when $\lambda_1 = \lambda_2 = \dots = \lambda_{m-1} > \lambda_m$ with $m > 2$. Letting $\alpha = \lambda_i^2(A)$ for any $i = 1, \dots, m-1$ and $\beta = \lambda_m^2(A)$, we have
\begin{align*}
	\underbrace{\frac{\beta}{m \alpha}}_{\text{U$_{\infty}$}} 
	< \underbrace{\frac{\beta}{\alpha(m-1) + \beta}}_{\text{NU}} 
	< \underbrace{\frac{\beta}{\alpha + \beta(m-1)}}_{\text{MR}_\infty} 
	\le \underbrace{\frac{1}{\lambda_{i_k}^2}\frac{\alpha \beta}{\alpha + \beta(m-1)}}_{\text{MR}} 
	< \underbrace{\frac{1}{m}}_{\text{U, MD}}.
\end{align*}
Thus, Strohmer and Vershynin's NU rule would actually be the worst rule to use, whereas U and MD are the best.
In this case 
$\const{}^2$ is closer to its upper bound ($\approx \beta$) so we would expect greedy rules to perform well.

\section{Approximate Greedy Rules}\label{sec:approx}

In many applications computing the exact MR or MD rule will be too inefficient, but we can always approximate it using a cheaper {\em approximate} greedy rule, as in the method of~\citet{eldar2011johnsonlindenstrauss}. In this section we consider methods that compute the greedy rules up to multiplicative or additive errors.

\subsection{Multiplicative Error}\label{subsec:multiplicative}

Suppose we have approximated the MR rule such that there is a multiplicative error in our selection of $i_k$,
\[
	|a_{i_k}^T x^k - b_{i_k} | \ge \max_i |a_i^Tx^k - b_i| (1 - \epsilon_k),
\]
for some $\epsilon_k \in [0,1)$. In this scenario, using the tight analysis for the MR rule, we show in Appendix~\ref{app:multiplicative} that 
\[
	\| x^{k+1} - x^* \|^2 \le \bigg ( 1 - \frac{(1 - \epsilon_k)^2 \const{}^2}{\| a_{i_k} \|^2} \bigg ) \| x^k - x^* \|^2.
\]
Similarly, if we approximate the MD  rule up to a multiplicative error,
\[
	\left | \frac{a_{i_k}^T x^k - b_{i_k}}{\|a_{i_k}\|} \right | \ge \max_i \left | \frac{a_i^Tx^k - b_i}{\|a_i\|} \right| (1 - \bar{\epsilon}_k),
\]
for some $\bar{\epsilon}_k \in [0,1)$, then we show in Appendix~\ref{app:multiplicative} that the following rate holds,
\[
	\| x^{k+1} - x^* \|^2 \le \bigg ( 1 - (1 - \bar{\epsilon}_k)^2 \constbar{}^2 \bigg ) \| x^k - x^* \|^2.
\]
These scenarios do not require the error to converge to $0$. However, if $\epsilon_k$ or $\bar{\epsilon}_k$ is large, then the convergence rate will be slow. 

\subsection{Additive Error}\label{subsec:additive}

Suppose we select $i_k$ using the MR rule up to additive error,
\[
	|a_{i_k}^T x^k - b_{i_k}|^2 \ge \max_i |a_i^Tx^k - b_i|^2 - \epsilon_k,
\]
or similarly for the MD rule,
\[
	\left |\frac{a_{i_k}^T x^k - b_{i_k}}{\|a_{i_k}\|} \right|^2 \ge \max_i \left |\frac{a_i^Tx^k - b_i}{\|a_i\|} \right |^2 - \bar{\epsilon}_k,
\]
for some $\epsilon_k \ge 0$ or $\bar{\epsilon}_k \ge 0$, respectively.
We show in Appendix~\ref{app:additive} that this results in the following convergence rates for the MR and MD rules with additive error (respectively),
\[
    \| x^{k+1} - x^* \|^2 \le \bigg ( 1 - \frac{\const{}^2}{\| a_{i_k} \|^2} \bigg ) \| x^k - x^* \|^2 + \frac{\epsilon_k}{\|a_{i_k}\|^2},
\]
and
\[
    \| x^{k+1} - x^* \|^2 \le \big ( 1 -   \constbar{}^2 \big ) \| x^k - x^* \|^2 + \bar{\epsilon_k}.
\]
With an additive error, we need the errors to go to 0 in order for the algorithm to converge; if it does go to $0$ fast enough, we obtain the same rate as if we were calculating the exact greedy rule. In the approximate greedy rule used by~\citet{eldar2011johnsonlindenstrauss}, there is unfortunately a constant additive error. To address this, they compare the approximate greedy selection to a randomly selected $i_k$ and take the one with the largest distance. This approach can be substantially faster when far from the solution, but may eventually revert to random selection. We give details comparing Eldar and Needell's rate to our above rate in Appendix~\ref{app:73}, but here we note that the above bounds will typically be much stronger.


\section{Systems of Linear Inequalities}\label{sec:inequalities}

Kaczmarz methods have been extended to systems of linear inequalities, 
\begin{equation}\label{ineqsystem}
	\begin{cases}
		a_i^Tx \le b_i & (i \in I_\le) \\
		a_i^Tx = b_i & (i \in I_=).
	\end{cases}
\end{equation}
where the disjoint index sets $I_\le$ and $I_=$ partition the set $\{ 1,2,\dots,m\}$~\citep{leventhal2010constraints}. 
In this setting the method takes the form
\begin{align*}
	x^{k+1} = x^k - \frac{\beta^k}{\| a_i \|^2} a_i,\quad
\text{with }\quad
	\beta^k  =
	\begin{cases}
		(a_i^T x^k - b_i)^+ & (i \in I_\le) \\
		~a_i^T x^k - b_i & (i \in I_=),
	\end{cases} 
\end{align*}
where $(\gamma)^+ = \max \{ \gamma,0 \}$. In Appendix~\ref{app:inequalities} we derive analogous greedy rules and convergence results for this case. The main difference in this setting is that the rates are in terms of the distance of $x^{k}$ to the feasible set $S$ of \eqref{ineqsystem},
\[
	d(x^{k},S) = \min_{z \in S} \| x^{k} - z \|_2 = \| x^k - P_{S}(x^k) \|_2,
\]
where $P_{S}(x)$ is the projection of $x$ onto $S$. This generalization is needed because with inequality constraints the different iterates $x^k$ may have different projections onto $S$.



\section{Multi-Step Analysis}\label{sec:multi}

All existing analyses of Kaczmarz methods consider convergence rates that depend on a {\em single} step (in the case of randomized/greedy selection rules) or a single cycle (in the cyclic case). In this section we derive the first tighter {\em multi}-step convergence rates for iterative Kaczmarz methods; we first consider the MR rule, and then we explore the potential of faster random selection rules. These new rates/rules depend on the orthogonality graph introduced in Section~\ref{sec:calc}, and thus in some sense they depend on the `angle' between rows. This dependence on the `angle' is similar to the classic convergence rate analyses of cyclic Kaczmarz algorithms, and is a property that is not captured by existing randomized/greedy analyses (which only depend on the row norms). 


\subsection{Multi-Step Maximum Residual Bound}\label{subsec:multi-MR}

If two rows $a_i$ and $a_j$ are orthogonal, then if the equality $a_{i}^Tx^k = b_{i}$ holds at iteration $x^k$ and we select $i_k = j$, then we know that $a_i^Tx^{k+1} = b_i$. More generally, updating $i_k$ makes equality $i_k$ satisfied but could make any equality $j$ unsatisfied where $a_j$ is not orthogonal to $a_{i_k}$.
Thus, after we have selected row $i_k$, equation $i_k$ will remain satisfied for all subsequent iterations until one of its neighbours is selected in the orthogonality graph. During these subsequent iterations, it cannot be selected by the MR rule since its residual is zero.

In Appendix~\ref{app:multi-MR}, we show 
how the structure of the orthogonality graph can be used to derive a worst-case bound on the \emph{sequence} of $\norm{a_{i_k}}$ values that appear in the tighter analysis of the MR rule~\eqref{maxrestight}. In particular, we show that the MR rule achieves a convergence rate of 
\begin{align*}
	\| x^k - x^* \|^2 \le O(1) \!\! \left(\max_{S(G)}\left\{\!\!\!\sqrt[\leftroot{-2}\uproot{2}|S(G)|]{\prod_{j \in S(G)} \left ( 1 - \frac{\const{}^2}{\|a_{j} \|^2} \! \right )} \right\}\right)^k \!\!\norm{x^0-x^*}^2,
\end{align*}
where  the maximum is taken over the geometric means of all the \emph{star subgraphs} $S(G)$ of the orthogonality graph with at least two nodes (these are the connected subgraphs that have a diameter of $1$ or $2$). Although this result is quite complex, even to state, there is a simple implication of it: if the values of $\norm{a_i}$ that are close to $\norm{A}_{\infty,2}$ are all more than two edges away from each other in the orthogonality graph, then the MR rule converges substantially faster than the worst-case MR$_\infty$ bound~\eqref{maxres} indicates. 

A multi-step analysis of coordinate descent with the Gauss-Southwell rule and exact coordinate optimization was recently considered by~\citet{nutini2015}. To derive this bound, they convert the problem to the same weighted graph construction we use in Appendix~\ref{app:multi-MR}. However, they were only able to derive a bound on this construction in the case of chain-structured graphs. Our result uses a generalization of their result to the case of general graphs, and indeed our result is \emph{tighter} than the bound that they conjectured would hold for general graphs. Since the graph construction in this work is the same as in their work, the result we use also gives the tightest known bound on coordinate descent with the Gauss-Southwell rule and exact coordinate optimization.


\subsection{Faster Randomized Kaczmarz Rules}\label{subsec:fasterrandom}

The orthogonality graph can also be used to design faster randomized algorithms. To do this, we use the same property as in the previous section: after we have selected $i_k$, equality $i_k$ will be satisfied on all subsequent iterations until we select one of its neighbours in the orthogonality graph. Based on this, we call a row $i$ `selectable' if $i$ has never been selected or if a neighbour of $i$ in the orthogonality graph has been selected since the last time $i$ was selected.\footnote{If we initialize with $x^0=0$, then instead of considering all nodes as initially selectable we can restrict to the nodes $i$ with $b_i \neq 0$ since otherwise we have $a_i^Tx^0 = b_i$ already.} We use the notation $s_i^k = 1$ to denote that row $i$ is `selectable' on iteration $k$, and otherwise we use $s_i^k = 0$ and say that $i$ is `not selectable' at iteration $k$. There is no reason to ever update a `not selectable' row, because by definition the equality is already satisfied. Based on this, we propose two simple randomized schemes:
\begin{enumerate}
\item {\bf Adaptive Uniform}: select $i_k$ uniformly from the selectable rows.
\item {\bf Adaptive Non-Uniform}: select $i_k$ proportional to $\norm{a_i}^2$ among the selectable rows.
\end{enumerate}
Let $A_k/\bar{A}_k$ denote the sub-matrix of $A/\bar{A}$ formed by concatenating the selectable rows on iteration $k$, and let $m_k$ denote the number of selectable rows. If we are given the set of selectable nodes at iteration $k$, then for adaptive uniform we obtain the bound
\[
	\mathbb{E} [\| x^{k+1} - x^* \|^2] \le \left ( 1 - \frac{\constobark{}^2}{m_k} \right ) \| x^k - x^* \|^2,
\]
while for adaptive non-uniform we obtain the bound
\[
	\mathbb{E}[\| x^{k+1} - x^* \|^2] \le \bigg ( 1 - \frac{\constok{}^2}{\| A_k \|^2_F} \bigg ) \| x^k - x^* \|^2.
\]
If we are not on the first iteration, then at least one node is not selectable and these are strictly faster than the previous bounds. The gain will be small if most nodes are selectable (which would be typical of dense orthogonality graphs), but the gain can be very large if only a few nodes are selectable (which would be typical of sparse orthogonality graphs).

\textbf{Theoretical Rate}: If we form a vector $s^k$ containing the values $s_i^k$, it's possible (at least theoretically) to compute the expected value of $s^k$ by viewing it as a Markov chain. In particular, $s^0$ is a vector of ones while $p(s^{k+1} | s^k)$ is equal to the normalized sum of all ways $s^{k+1}$ could be the set of selectable nodes given the selectable nodes $s^k$ and the orthogonality graph (most $p(s^{k+1} | s^k)$ values will be zero). Given this definition, we can express the probability of a particular $s^k$ recursively using the Chapman-Kolmogorov equations,
\[
	p(s^{k+1}) = \sum_{s^k} p(s^{k+1} | s^k) p(s^k).
\]
If we are interested in the probability that a particular $s_i^k = 1$, we can sum $p(s^k)$ over values $s^k$ compatible with this event. Unfortunately, deriving tighter bounds using these probabilities appears to be highly non-trivial.

\textbf{Practical Issues}:
In order for the adaptive methods to be efficient, we must be able to efficiently form the orthogonality graph and update the set of selectable nodes. If each node has at most $g$ neighbours in the orthogonality graph, then the cost of updating the set of selectable nodes and then sampling from the set of selectable nodes is $O(g\log(m))$  (we give details in Appendix~\ref{app:fasterrandom}). In order for this to not increase the iteration cost compared to the NU method, we only require the very-reasonable assumption that $g\log(m)  = O(n + \log(m))$. In many applications where orthogonality is present, $g$ will be far smaller than this.

However, forming the orthogonality graph at the start may be prohibitive: it would cost $O(m^2n)$ in the worst case to test pairwise orthogonality of all nodes. In the sparse case where each column has at most $c$ non-zeros, we can find an approximation to the orthogonality graph in $O(c^2n)$ by assuming that all rows which share a non-zero are non-orthogonal. Alternately, in many applications the orthogonality graph is easily derived from the structure of the problem. For example, in graph-based semi-supervised learning where the graph is constructed based on the $k$-nearest neighbours, the orthogonality graph will simply be the given $k$-nearest neighbour graph as these correspond to the columns that will be mutually non-zero in $A$.


\section{Experiments}\label{sec:numerics}

\begin{figure*}
\centering
\includegraphics[width=.49\textwidth]{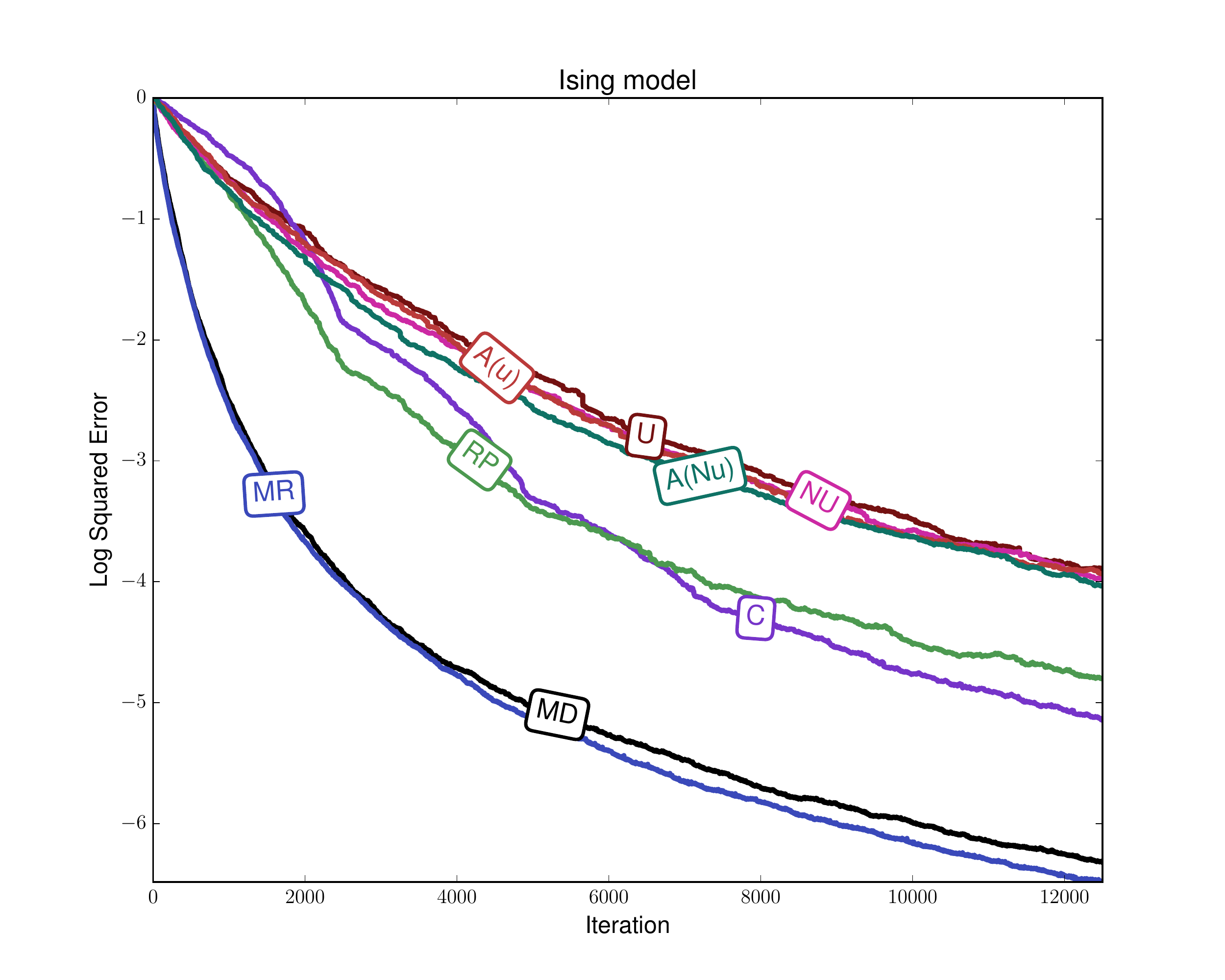}
\includegraphics[width=.49\textwidth]{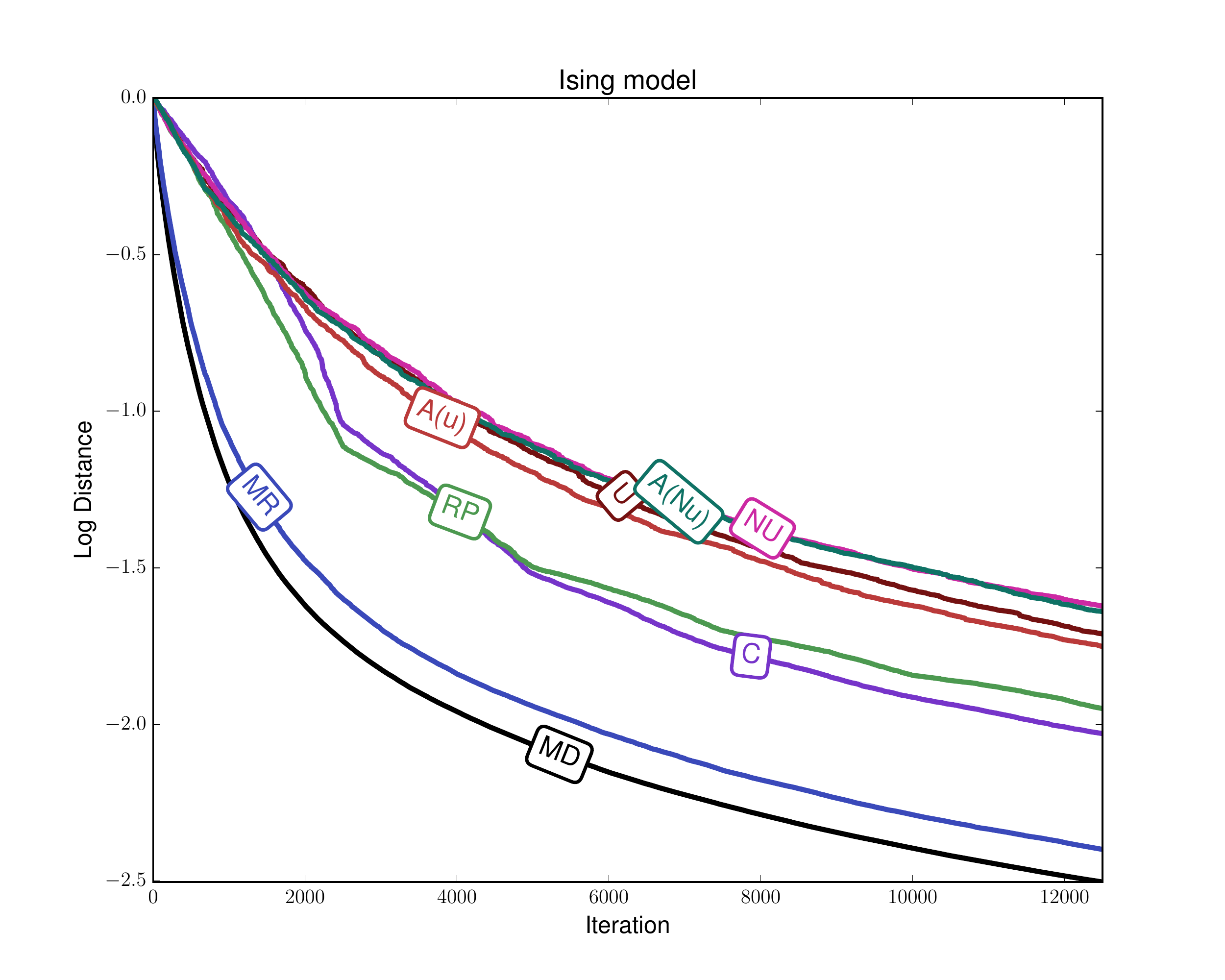}\\ 
\includegraphics[width=.49\textwidth]{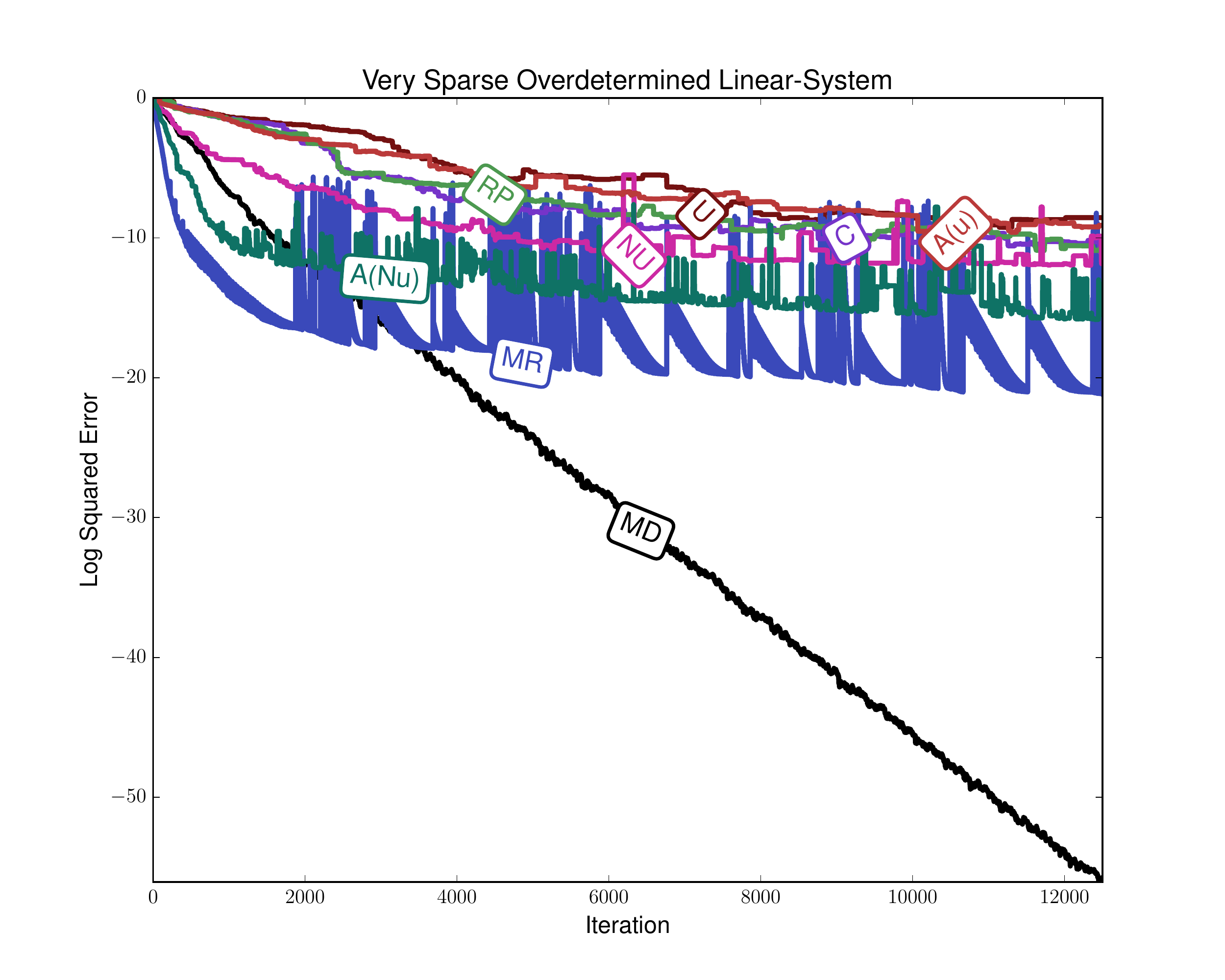}
\includegraphics[width=.49\textwidth]{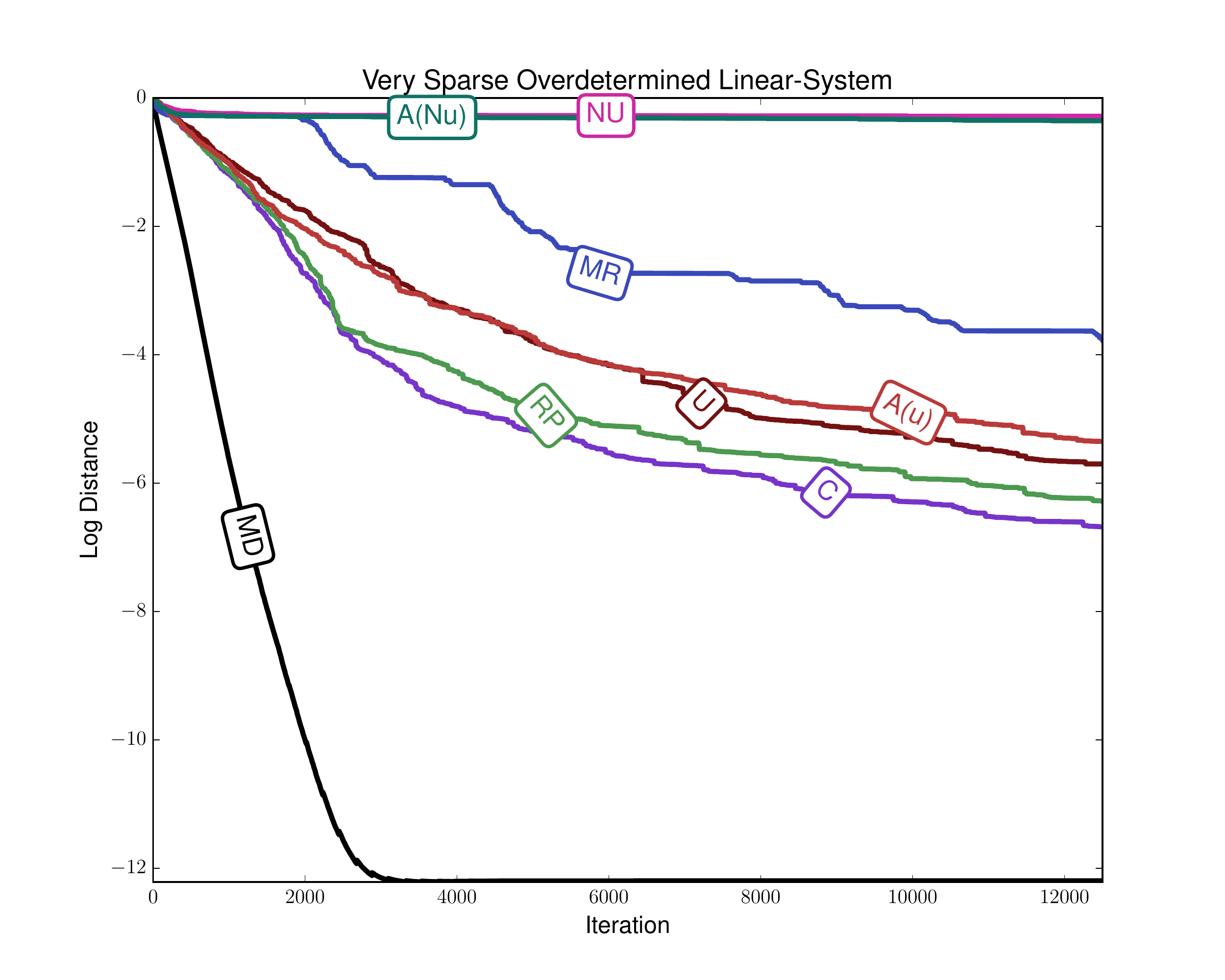}\\ 
\includegraphics[width=.49\textwidth]{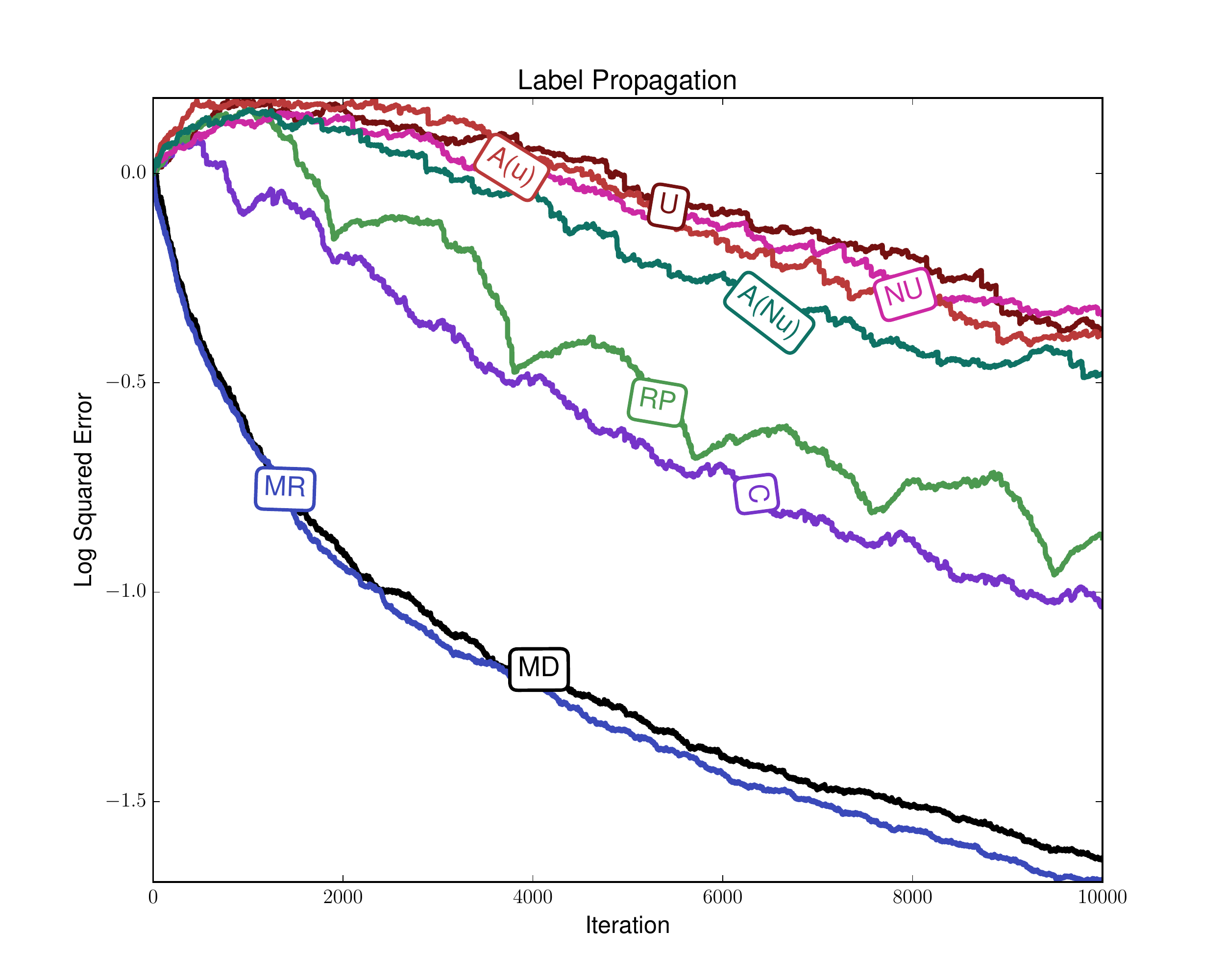}  
\includegraphics[width=.49\textwidth]{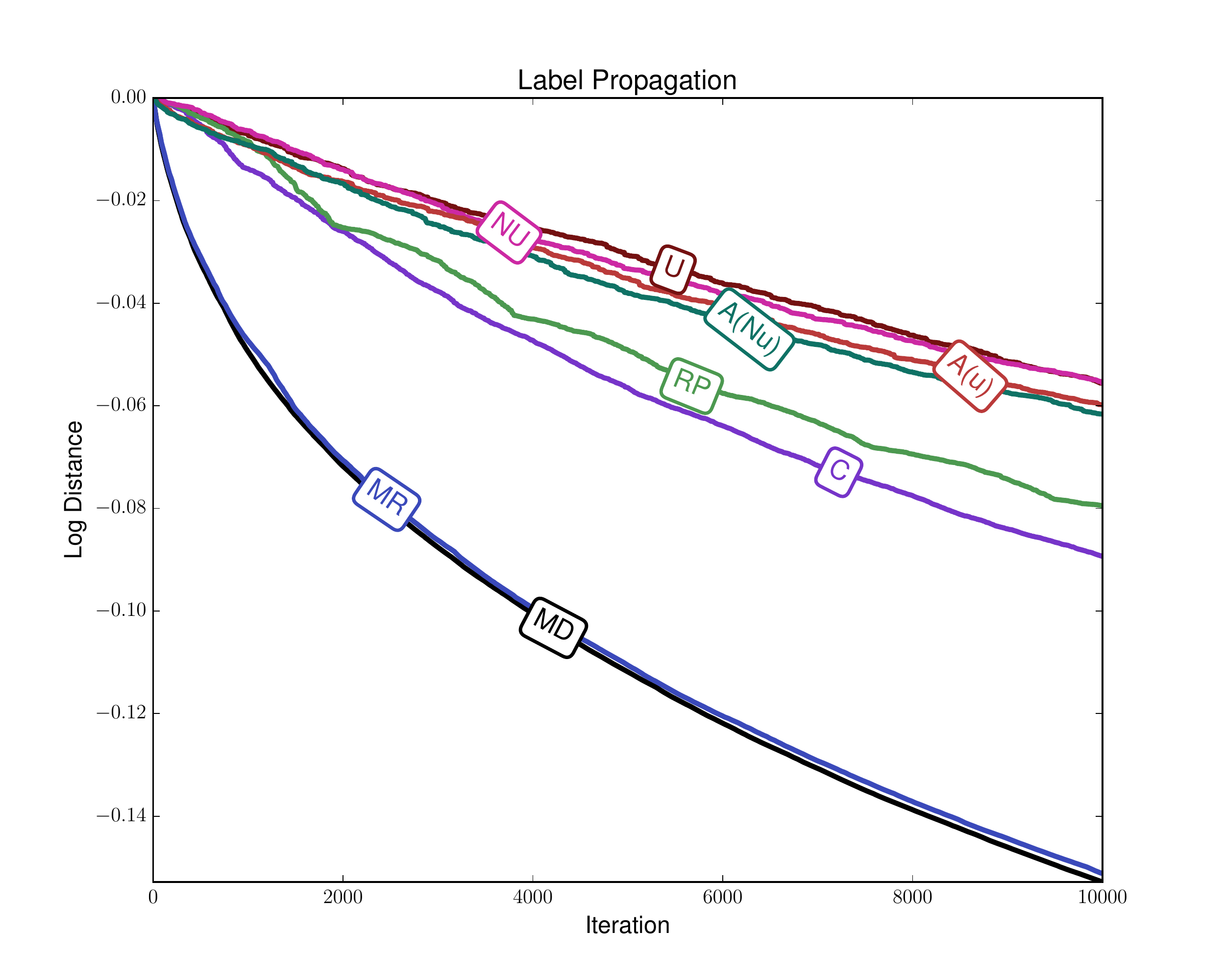} 
\caption{Comparison of Kaczmarz selection rules for squared error (left) and distance to solution (right).}
\label{fig:h1}
\end{figure*}

\citet{eldar2011johnsonlindenstrauss} have already shown that approximate greedy rules can outperform randomized rules for  dense problems. Thus, in our experiments we focus on comparing the effectiveness of different rules on very sparse problems where our max-heap strategy allows us to efficiently compute the exact greedy rules. The first problem we consider is generating a dataset $A$ with a $50$ by $50$  lattice-structured dependency (giving $n = 2500$). The corresponding $A$ has the following non-zero elements: the diagonal elements $A_{i,i}$, the upper/lower diagonal elements $A_{i,i+1}$ and $A_{i+1,i}$ when $i$ is not a multiple of $50$ (horizontal edges), and the diagonal bands $A_{i,i+50}$ and $A_{i+50,i}$ (vertical edges). We generate these non-zero elements from a $\mathcal{N}(0,1)$ distribution
and generate the target vector $b = Az$ using  $z \sim \mathcal{N}(0,I)$.
Each row in this problem has at most four neighbours, and this type of sparsity structure is typical of spatial Gaussian graphical models and linear systems that arise from discretizing two-dimensional partial differential equations.


The second problem we consider is solving an overdetermined consistent linear system with a very sparse $A$ of size $2500 \times 1000$. We generate each row of $A$ independently such that there are ${\log(m)}/{2m}$ non-zero entries per row drawn from a uniform distribution between 0 and 1. To explore how having different row norms affects  the performance of the selection rules, we randomly multiply one out of every 11 rows by a factor of 10,000.

For the third problem, we solve a label propagation problem for semi-supervised learning in the `two moons' dataset \citep{zhou2004learning}. From this dataset, we generate 2000 samples and randomly label 100 points in the data. We then connect each node to its 5 nearest neighbours. Constructing a data set with such a high sparsity level is typical of graph-based methods for semi-supervised learning. We use a variant of the quadratic labelling criterion of \citet{bengio2006},
\[
	\min_{y_i | i \not \in S} \quad \frac{1}{2} \sum_{i=1}^n \sum_{j=1}^n w_{ij} (y_i - y_j)^2,
\]
where $y$ is our label vector (each $y_i$ can take one of 2 values), $S$ is the set of labels that we do know and $w_{ij} \ge 0$ are the weights assigned to each $y_i$ describing  how strongly we want the label $y_i$ and $y_j$ to be similar. We can express this quadratic problem as a linear system that is consistent by construction (see Appendix~\ref{app:numerics}), and hence apply Kaczmarz methods. As we labelled $100$ points in our data, the resulting linear system has a matrix of size $1900 \times 1900$ while the number of neighbours $g$ in the orthogonality graph is at most $5$.

In Figure \ref{fig:h1} we compare the normalized squared error and distance against the iteration number for 8 different selection rules: cyclic (C), random permutation (RP - where we change the cycle order after each pass through the rows), uniform random (U), adaptive uniform random (A(u)), non-uniform random NU, adaptive non-uniform random (A(Nu)), maximum residual (MR), and maximum distance (MD).

In experiments 1 and 3, MR performs similarly to MD and both outperform all other selection rules. For experiment 2, the MD rule outperforms all other selection rules in terms of distance to the solution although MR performs better on the early iterations in terms of squared error. In Appendix~\ref{app:numerics} we explore a `hybrid' method on the overdetermined linear system problem that does well on both measures. In Appendix~\ref{app:numerics}, we also plot the performance in terms of runtime.


The new randomized A(u) method did not give significantly better performance than the existing U method on any dataset. This agrees with our bounds which show that the gain of this strategy is modest. In contrast, the new randomized A(Nu) method performed much better than the existing NU method on the over-determined linear system in terms of squared error. This again agrees with our bounds which suggest that the A(Nu) method has the most to gain when the row norms are very different. Interestingly, in most experiments we found that \emph{cyclic} selection worked better than any of the randomized algorithms. However, cyclic methods were clearly beaten by greedy methods.



\section{Discussion}

In this work, we have proven faster convergence rate bounds for a variety of row-selection rules in the context of Kaczmarz methods for solving linear systems. We have also provided new randomized selection rules that make use of orthogonality in the data in order to achieve better theoretical and practical performance.  While
we have focused on the case of non-accelerated and single-variable variants of the Kaczmarz algorithm,  we expect that all of our conclusions also hold for accelerated Kaczmarz and block Kaczmarz methods \citep{needell2012, lee2013, liu2014accelerated, gower2015, oswald2015}.




\subsubsection*{Acknowledgements}

This research was supported by the Natural Sciences and Engineering Research Council of Canada (NSERC RGPIN-06068-2015), and 
Julie Nutini is funded by an NSERC Canada Graduate Scholarship.

\appendix
\renewcommand{\thesubsection}{\Alph{subsection}}
\section*{Appendix}

\section{Efficient Calculations for Sparse $A$}\label{app:calc} 

To compute the MR rule efficiently for sparse $A$, we need to store and update the residuals $r_i = (a_i^Tx^k - b_i)$ for all $i$. If we initialize with $x^0 = 0$, then the initial values of the residuals are simply the corresponding $b_i$ values.
Given the initial residuals, we can construct a max-heap structure on these residuals in $O(m)$ time.
 The max-heap structure lets us compute the MR rule in $O(1)$ time. After an iteration of the Kaczmarz method, we can update the max-heap
efficiently as follows:

{\bf For each $j$ where $x_j^{k+1} \neq x_j^k$:}
\vspace{-0.5em}
\begin{itemize}
		\item {\bf For each $i$ with $a_{ij} \not = 0$:}
\begin{itemize}
\item Update $r_i$ using $r_i \leftarrow r_i - a_{ij}x^k_j + a_{ij}x^{k+1}_j.$
\item Update max-heap using the new value of $|r_i|$.
\end{itemize}
\end{itemize}
The cost of each update to an $r_i$ is $O(1)$ and the cost of each heap update is $O(\log m)$. If each row of $A$ has at most $r$ non-zeroes and each column has at most $c$ non-zeroes, then the outer loop is run at most $r$ times while the inner loop is run at most $c$ times for each outer loop iteration. Thus, in the worst case the total cost is $O(cr \log m)$, although it might be much faster if we have particularly sparse rows or columns. Thus, if $c$ and $r$ are sufficiently small, the MR rule is not much more expensive than non-uniform random selection which costs $O(\log m)$. For the MD rule, the cost is the same except there is  an extra one-time cost to pre-compute the row norms $\|a_i\|$. 

Now consider the case where $A$ may be dense but each row is orthogonal to all but at most $g$ other rows. In this setting it would be too slow to implement the above update of the residuals, since the cost would be $O(mn\log(m))$. In this setting, it makes more sense to use the following alternative approach to update the max-heap after we've updated row $i_k$:

{\bf For each $i$ that is a neighbour of $i_k$ in the orthogonality graph:}
\begin{itemize}
\item Compute the residual $r_i = a_i^Tx^k - b_i$.
\item Update max-heap using the new value of $|r_i|$.
\end{itemize}
We can find the set of neighbours for each node in constant time by keeping a list of each node's neighbours. This loop would run at most $g$ times and the cost of each iteration would be $O(n)$ to update the residual and $O(\log m)$ to update the heap. Thus, the cost to track the residuals using this alternative approach would be $O(gn + g\log(m))$ or the faster $O(gr + g\log(m))$ if each row has at most $r$ non-zeros.

\section{Randomized and Maximum Residual}\label{app:UandMR}

In this section, we provide details of the convergence rate derivations for the non-uniform and maximum residual (MR) selection rules. All the convergence rates we discuss use the following relationship,
\begin{align*}
	\| x^{k+1} - x^* \|^2
	& = \| x^k - x^* \|^2 - \|x^{k+1} - x^k \|^2 \\
	& = \| x^k - x^* \|^2 - \left \| \frac{(b_i - a_i^T x^k)}{\|a_i\|^2} \cdot a_i \right \|^2 \\
	&~= \| x^k - x^* \|^2 - \frac{\left (a_i^Tx^k - b_i \right)^2}{\| a_i \|^2} \numberthis \label{eq:eq6},
\end{align*}
which is equation \eqref{derivationstart} in the main paper.

\subsection*{Non-Uniform}

We first review how the steps discussed by \citet{vishnoi2012laplacian} that can be used to derive the convergence rate bound of \citet{strohmer2009} for non-uniform random selection when row $i$ is chosen according to the probability distribution determined by $\|a_i \|/ \|A \|_F$. Taking the expectation of \eqref{eq:eq6} with respect to $i$, we have
\begin{align*}
	\mathbb{E} [\| x^{k+1} - x^* \|^2]
	&~= \| x^k - x^* \|^2 - \mathbb{E} \bigg [  \frac{(a_i^T x^k - b_i)^2}{\| a_i \|^2} \bigg ]	\\
	&~= \| x^k - x^* \|^2 - \sum_{i = 1}^m \frac{\| a_i \|^2}{\| A \|^2_F} \frac{(a_i^\top (x^k - x^*))^2}{\| a_i \|^2} \\
	&~= \| x^k - x^* \|^2 - \frac{\| A (x^k - x^*) \|^2}{\| A \|_F^2} \\
	&~\le \left ( 1 - \frac{\consto{}^2}{\| A \|_F^2} \right ) \| x^k - x^* \|^2, \numberthis \label{nonuniformrate}
\end{align*}
where $\sigma(A,2)$ is the Hoffman~\citeyearpar{hoffman1952} constant, which can be defined as the largest value such that for any $x$ that is not a solution to the linear system we have
\begin{equation}\label{eq:hoffman2}
\sigma(A,2)\norm{x - x^*} \leq \norm{A(x - x^*)},
\end{equation}
where $x^*$ is the projection of $x$ onto the set of solutions $S$. In other words, we can write it as
\[
	\consto{} := \inf_{x \not\in S}\frac{\| A(x - x^*) \|}{\| x - x^* \|}.
\]
 \citet{strohmer2009} consider the special case where $A$ has independent columns, and this result yields their rate in this special case since under this assumption  $\sigma(A,2)$ is given by the $n$th singular value of $A$. For general matrices, $\sigma(A,2)$ is given by the smallest non-zero singular value of $A$.

\subsection*{Maximum Residual}

We use a similar analysis to prove a convergence rate bound for the MR rule,
\begin{equation}\label{eqapp:maxResRule}
	i_k = \argmax{i} ~|a_i^T x^k - b_i|.
\end{equation}
Assuming that $i$ is selected according to \eqref{eqapp:maxResRule}, then starting from \eqref{eq:eq6} we have
\begin{align*}
    \| x^{k+1} - x^* \|^2
        &~= \| x^k - x^* \|^2 -  \frac{\max_i(a_i^T x^k - b_i)^2}{\| a_i \|^2}	\\
        &\le \| x^k - x^* \|^2 - \frac{1}{\| A \|^2_{\infty,2}} \max_i ~(a_i^T (x^k - x^*))^2 \\
	&= \| x^k - x^* \|^2 - \frac{\| A(x^k - x^*) \|^2_\infty}{\| A \|^2_{\infty,2}}  \\
	&\le \bigg ( 1 - \frac{\const{}^2}{\| A \|^2_{\infty,2}} \bigg ) \| x^k - x^* \|^2 \numberthis \label{maxResRate},
\end{align*}
where $\| A \|_{\infty,2}^2 := \max_i \{ \| a_i \|^2 \}$ and $\const{}$ is the largest value such that
\begin{equation}\label{eq:hoffmaninf}
\const{}\norm{x-x^*} \leq \norm{A(x-x^*)}_\infty,
\end{equation}
or equivalently
\[
\const{} := \inf_{x \not\in S}\frac{\norm{A(x - x^*)}_\infty}{\norm{x-x^*}}.
\]
The existence of such a Hoffman-like constant follows from the existence of the Hoffman constant and the equivalence between norms. Applying the norm equivalence $\| \cdot \|_\infty \ge \frac{1}{\sqrt{m}} \| \cdot \|$ to equation \eqref{eq:hoffman2} we have
\[
	\sigma(A,2)\norm{x - x^*} \leq \norm{A(x - x^*)} \leq \sqrt{m} \norm{A(x - x^*)}_\infty,
\]
which implies that $\sigma(A,2)/\sqrt{m} \leq \sigma(A,\infty)$.
Similarly, applying $\| \cdot \|_\infty \le \| \cdot \|$ to \eqref{eq:hoffmaninf} we have
\[
	\const{}\norm{x-x^*} \leq \norm{A(x-x^*)}_\infty \le \norm{A(x-x^*)},
\]
which implies that $\sigma(A,\infty)$ cannot be larger than $\sigma(A,2)$.
Thus, $\const{}$ satisfies the  relationship
\begin{equation}\label{sigma:Ainf}
	\frac{\consto{}}{\sqrt{m}} \le \const{} \le \consto{}.
\end{equation}

\section{Tighter Uniform and MR Analysis}\label{app:tight}

To avoid using the inequality $\|a_i \| \le \|A \|_{\infty,2}$ for all $i$, we want to `absorb' the individual row norms into the bound. We start with uniform selection.

\subsection*{Uniform}
Consider the diagonal matrix $D = \diag{\|a_1\|^2, \| a_2 \|^2, \dots, \| a_m \|^2}$. By taking the expectation of \eqref{eq:eq6}, we have
\begin{align*}
	\mathbb{E} [\| x^{k+1} - x^* \|^2] 
	&~= \| x^k - x^* \|^2 - \mathbb{E} \left [ \frac{\left (a_i^Tx^k - b_i \right)^2}{\| a_i \|^2} \right] \\
	&= \| x^k - x^* \|^2 - \sum_{i=1}^m \frac{1}{m} \frac{(a_i^T x^k - b_i)^2}{\| a_i \|^2}  \\
	&= \| x^k - x^* \|^2 - \frac{1}{m}\sum_{i=1}^m  \left (\left [\frac{a_i}{\| a_i \|} \right ]^T (x^k - x^*) \right)^2\\
	&= \| x^k - x^* \|^2 - \frac{\| D^{-1}A(x^k - x^*) \|^2}{m}  \\
	&\le \bigg ( 1  - \frac{\constobar{}^2}{m} \bigg ) \| x^k - x^* \|^2, \numberthis \label{eqapp:uniformtight}
\end{align*} 
where recall that $\bar{A} = D^{-1}A$ and we've used that $Ax = b$ and $D^{-1}Ax = D^{-1}b$ have the same solution set.

\subsection*{Maximum Residual}

For the tighter analysis of the MR rule we do not want to alter the selection rule. Thus, we first evaluate the MR rule and then divide by the corresponding $\| a_{i_{k}} \|^2$ for the selected $i_k$ at iteration $k$. Starting from \eqref{eq:eq6}, this gives us

\begin{align*}
    \| x^{k+1} - x^* \|^2
        &= \| x^k - x^* \|^2 - \frac{1}{\| a_{i_{k}} \|^2} \max_i (a_i^T(x^k - x^*))^2	\\
	&= \| x^k - x^* \|^2 - \frac{\| A(x^k - x^*) \|^2_\infty}{\| a_{i_{k}} \|^2} \\
	&\le \bigg ( 1 - \frac{\const{}^2}{\| a_{i_{k}} \|^2} \bigg ) \| x^k - x^* \|^2. \numberthis \label{tightmaxres}
\end{align*}
Applying this recursively over all $k$ iterations yields the  rate
\begin{equation}\label{rate:tightmaxres}
	\| x^k - x^* \|^2 \le \prod_{j = 1}^k\bigg ( 1 - \frac{\const{}^2}{\| a_{i_{j}} \|^2} \bigg ) \| x^0 - x^* \|^2.
\end{equation}

\section{Maximum Distance Rule}\label{app:maxDist}

If we can only perform one iteration of the Kaczmarz method, the {\em optimal} rule with respect to iterate progress is the maximum distance (MD) rule,
\begin{equation}\label{eqapp:maxdist}
	i_k = \argmax{i} \left | \frac{a_i^T x^k - b_i}{\|a_i \|} \right |.
\end{equation}

Starting again from \eqref{eq:eq6} and using $D$ as defined in the tight analysis for the U rule, we have
\begin{align*}
	\| x^{k+1} - x^* \|^2
	&= \| x^k - x^* \|^2 - \max_i  \bigg (\frac{a_i^T x^k -  b_i}{\|a_i\|} \bigg )^2 \\
	&= \| x^k - x^* \|^2 - \max_i  \bigg ( \left [\frac{a_i}{\|a_i\|}\right]^T (x^k -  x^*) \bigg )^2 \\
	&= \| x^k - x^* \|^2 -   \| D^{-1}A (x^k - x^*)  \|^2_\infty \\
	&\le \big ( 1 -   \constbar{}^2 \big ) \| x^k - x^* \|^2. \numberthis \label{maxDistRate}
\end{align*}
We now show that
\begin{equation}
\label{eq:distBound}
	\max \left \{ \frac{\constobar{}}{\sqrt{m}}, \frac{\consto{}}{\|A\|_F}, \frac{\const{}}{\|A\|_{\infty,2}} \right \}  \le \constbar{} \le \constobar{}.
\end{equation}
To derive the upper bound on $\constbar{}$, and to derive the lower bound in terms of $\constobar{}$, we can use norm equivalence arguments as we did for $\const{}$. This yields
\[
	\frac{\sigma(\bar{A},2)}{\sqrt{m}} \le \constbar{} \le \constobar{}.
\]

The last argument in the maximum in~\eqref{eq:distBound}, corresponding to the MR$_\infty$ rate, holds because \\$\|A \|_{\infty,2}~\ge~\|a_i\|$ for all $i$ so we have
\begin{align*}
\frac{\sigma(A,\infty)}{\norm{A}_{\infty,2}}\norm{x-x^*} 
&\leq \frac{\norm{A(x-x^*)}_\infty}{\norm{A}_{\infty,2}} \\
&= \max_i \left\{ \frac{|a_i^T(x-x^*)|}{\norm{A}_{\infty,2}}\right\} \\
&\le \max_i \left\{ \frac{|a_i^T(x-x^*)|}{\norm{a_i}}\right\} \\
&= \norm{\bar{A}(x-x^*)}_\infty.
\end{align*}

For the second argument in the maximum in~\eqref{eq:distBound}, the NU rate, we have
\begin{align*}
\frac{\sigma(A,2)^2}{\norm{A}_F^2}\norm{x-x^*}^2 
&\leq \frac{\norm{A(x-x^*)}^2}{\norm{A}_F^2} \\
&= \frac{\sum_{i}(a_i^T(x-x^*))^2}{\sum_i \norm{a_i}^2} \\
&\leq \max_i \left\{ \frac{(a_i^T(x-x^*))^2}{\norm{a_i}^2}\right\} \\
&= \norm{\bar{A}(x-x^*)}_\infty.
\end{align*}
The second inequality 
is true by noting that it is equivalent to the inequality
\[
1 \leq \max_i \left\{\frac{(a_i^T(x-x^*)^2/\sum_j (a_j^T(x-x^*))^2}{\norm{a_i}^2/\sum_j \norm{a_j}^2}\right\},
\]
and this true because the maximum ratio between two probability mass functions must be at least 1,
\[
	1 \le \max_i \frac{p_i / \sum_j p_j}{q_i / \sum_j q_j}, \quad \text{with all } p_i \ge 0, q_i \ge 0.
\]

Finally, we note that the MD rule  obtains the tightest bound in terms of performing one step. This follows from~\eqref{eq:eq6},
\[
\| x^{k+1} - x^* \|^2 = \| x^k - x^* \|^2 - \|x^{k+1} - x^k \|^2 = \| x^k - x^* \|^2 - \frac{\left (a_i^Tx^k - b_i \right)^2}{\| a_i \|^2},
\]
and noting that the MD rule maximizes $\|x^{k+1} - x^k \|$ and thus it maximizes how much smaller $\norm{x^{k+1}-x^*}$ is than  $\norm{x^k - x^*}$.

\section{Kaczmarz and Coordinate Descent}\label{app:cd}

Consider the Kaczmarz update:
\[
	x^{k+1} = x^k - \frac{(a_i^T x^k - b_i)}{\| a_i \|^2} a_i.
\]
This update is equivalent to one step of coordinate descent (CD) with step length $1/\|a_i\|^2$ applied to the dual problem,
\begin{equation}\label{eqapp:dual}
	\min_y \frac{1}{2} \| A^Ty \|^2 - b^Ty,
\end{equation}
see~\citet{wright2015}. 
%
%
Using the primal-dual relationship $A^Ty = x$, we can show the relationship between the greedy Kaczmarz selection rules and applying greedy coordinate descent rules to this dual problem. Consider the gradient of the dual problem,
\[
	\nabla f(y) = AA^Ty - b.
\]
The Gauss-Southwell (GS) rule for CD on the dual problem is equivalent to the MR rule for Kaczmarz on the primal problem since
\[
	i_k = \underbrace{\argmax{i} |\nabla_i f(y^k) |}_{\text{Gauss-Southwell rule}} = \argmax{i}  |a_i^T(A^T y^k) - b_i| = \underbrace{\argmax{i} |a_i^Tx^k - b_i|}_{\text{Maximum residual rule}}
\]
where $a_i^T$ is the $i$th row of $A$. 
Similarly, the Gauss-Southwell-Lipschitz (GSL) rule \citep{nutini2015} applied to the dual is equivalent to applying a Kaczmarz iteration with  the MD rule,
\[
	i_k = \underbrace{\argmax{i} \frac{|\nabla_i f(y^k) |}{\sqrt{L_i}}}_{\text{Gauss-Southwell-Lipschitz rule}} = \argmax{i}  \frac{|a_i^T(A^Ty^k) - b_i|}{\|a_i\|} = \underbrace{\argmax{i} \left | \frac{a_i^Tx^k - b_i}{\|a_i\|} \right |}_\text{Maximum distance rule},
\]
as the Lipschitz constants for the dual problem are $L_i = \|a_i \|^2$.

\begin{figure*}[h]
\centering
\includegraphics[width=.49\textwidth]{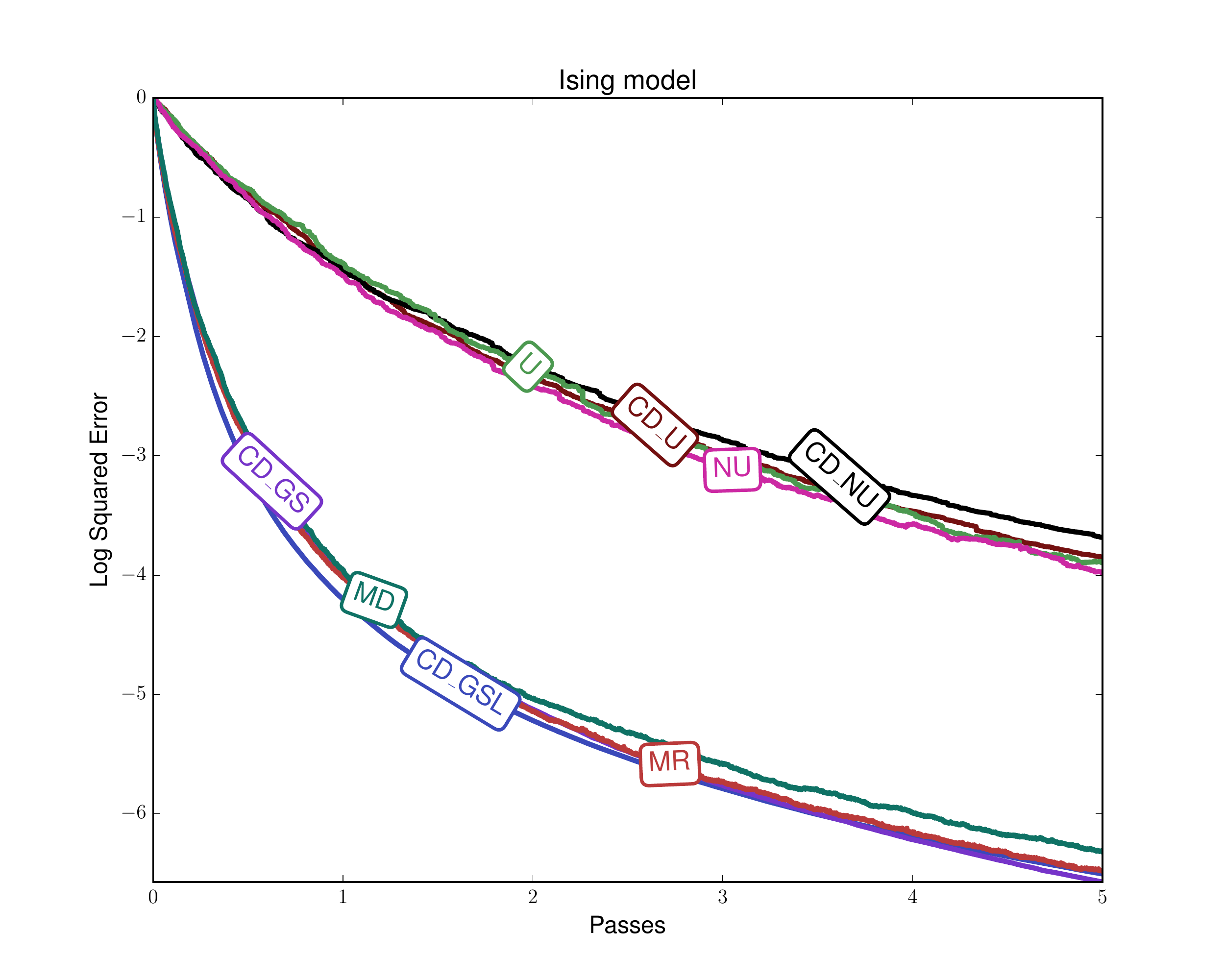}
\includegraphics[width=.49\textwidth]{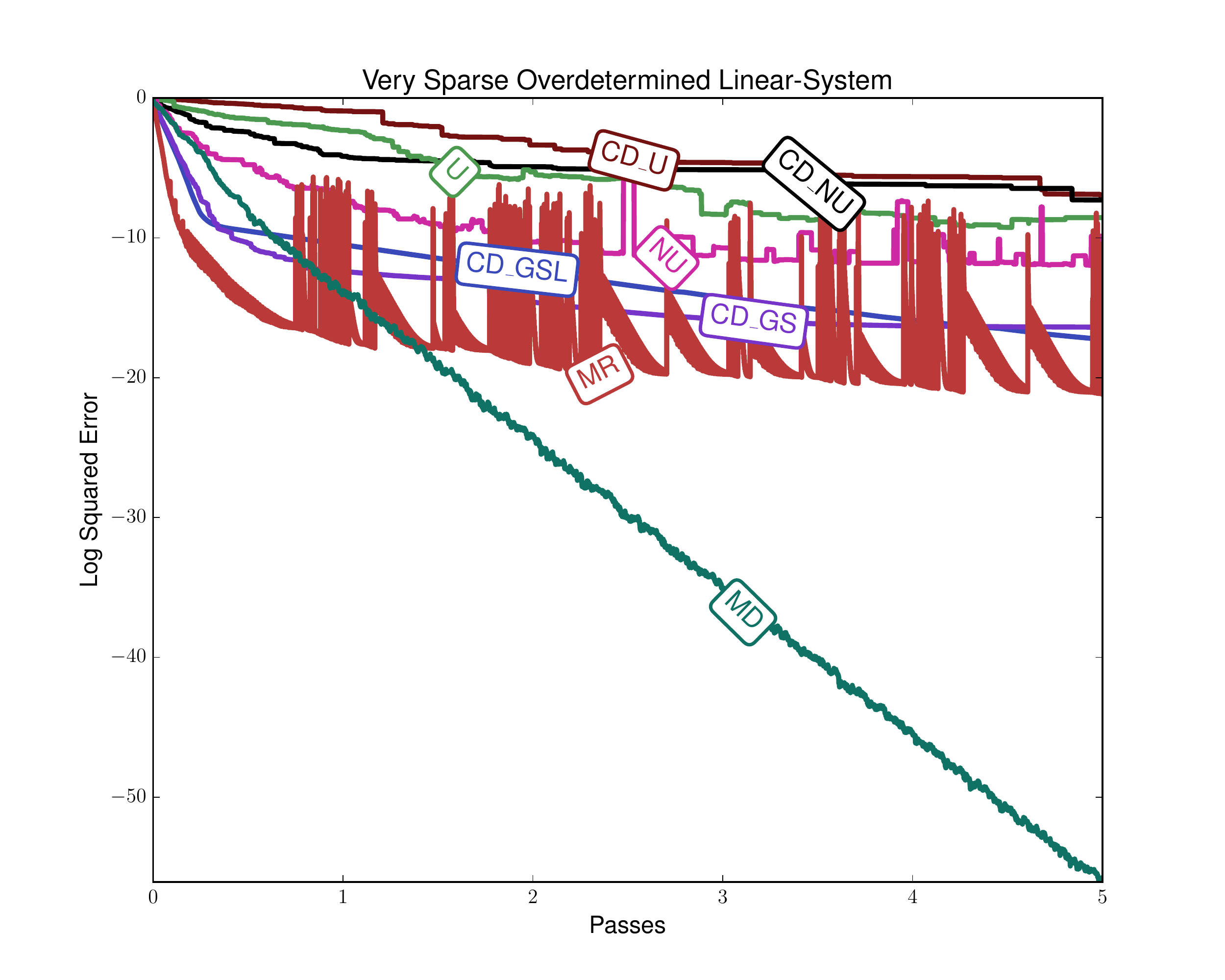}\\
\includegraphics[width=.49\textwidth]{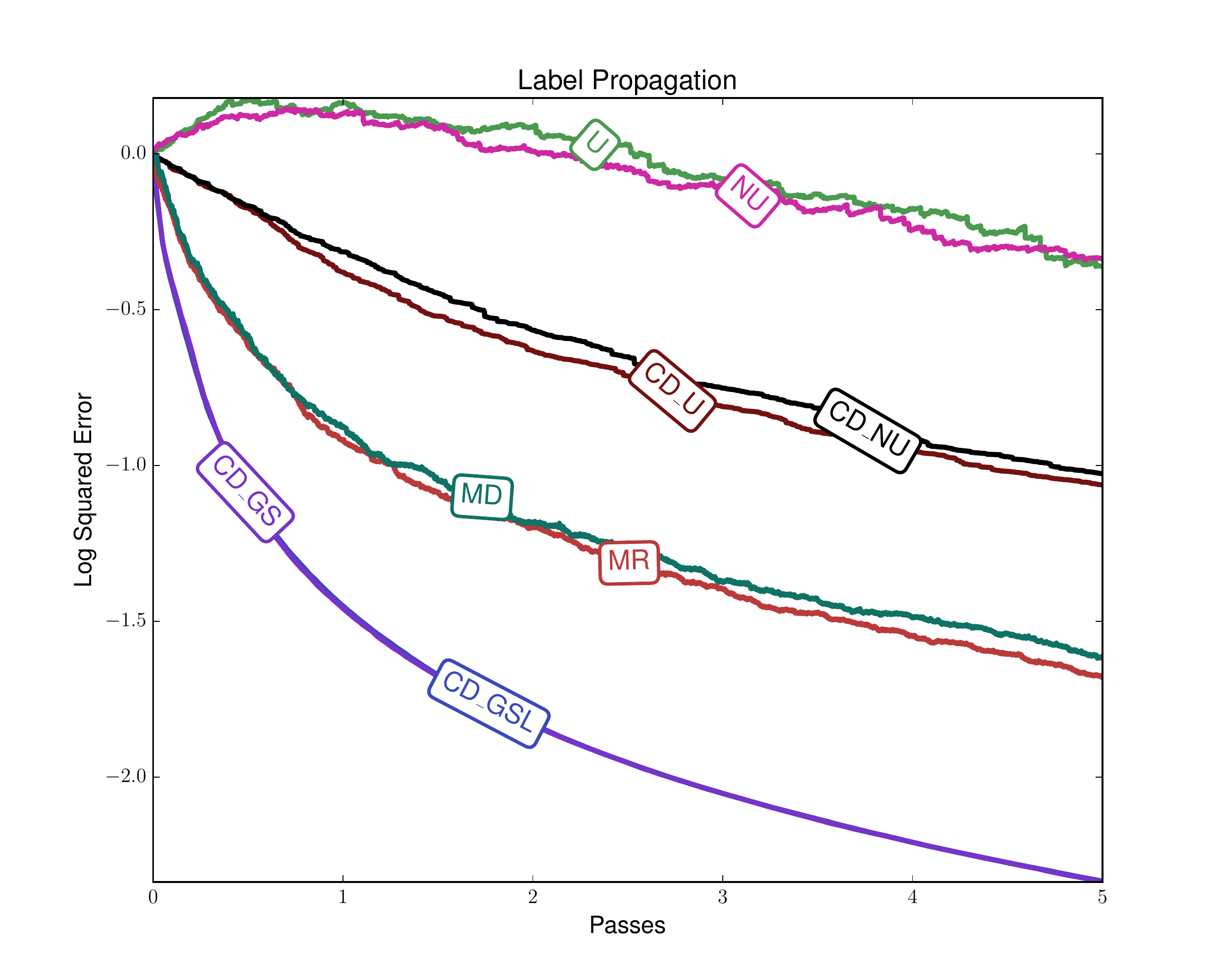}
\caption{Comparison of Kaczmarz and Coordinate Descent.}
\label{fig:KCD}
\end{figure*}
\bigskip
Figure \ref{fig:KCD} shows the results of running Kaczmarz compared to using CD (on the least-squares primal problem) for our 3 datasets from Section \ref{sec:numerics} of the main paper. In this figure we measure the performance in terms of the number of ``effective passes" through the data (one ``effective" pass would be the number of iterations needed for the cyclic variant of the algorithm to visit the entire dataset). In the first experiment Kaczmarz and CD methods perform similarly, while Kaczmarz methods work better in the second experiment and CD methods work better in the third experiment.

\section{Example: Diagonal $A$}\label{app:diagonal}

Consider a square diagonal matrix $A$ with $a_{ii} > 0$ for all $i$. In this case, the diagonal entries are the eigenvalues $\lambda_i$ of the $A$ and $\consto{} = \lambda_{\min}$. We give the convergence rate constants for such a diagonal $A$  in Table \ref{apptable:1}, and in this section we show how to arrive at these rates.
\renewcommand{\arraystretch}{2}
\begin{table}[h!]\
\centering
\caption{Convergence Rate Constants for Diagonal A}\label{apptable:1}
\begin{tabular}{| l | c | c |}
	\hline
	Rule & Rate &  Diagonal $A$\\ \hline
	U$_\infty$ & $\displaystyle \left (1 - \frac{\consto{}^2}{m \| A \|_{\infty,2}^2} \right)$ &  $\displaystyle \left (1 - \frac{\lambda_{\min}^2}{m \lambda_{\max}^2} \right )$ \\ \hline
	U & $\displaystyle \left (1 - \frac{\constobar{}^2}{m} \right)$ &  $\displaystyle \left (1 - \frac{1}{m} \right )$ \\ \hline
	NU & $\displaystyle \left (1 - \frac{\consto{}^2}{\| A \|_F^2}\right)$ & $\displaystyle \left (1 - \frac{\lambda_{\min}^2}{\sum_i \lambda_{i}^2} \right )$ \\ \hline
	MR$_\infty$ & $\displaystyle \left (1 - \frac{\const{}^2}{\| A \|_{\infty,2}^2}\right)$  & $\left ( 1 -  \displaystyle  \frac{1}{\lambda_1^2} \left [ \sum_i\frac{1}{\lambda_i^2} \right ]^{-1} \right )$ \\ \hline
	MR & $\displaystyle \left (1 - \frac{\const{}^2}{\| a_{i_k} \|^2}\right)$ & $\left ( 1 -  \displaystyle  \frac{1}{\lambda_{i_k}^2} \left [ \sum_i\frac{1}{\lambda_i^2} \right ]^{-1} \right )$ \\ \hline
	MD & $\displaystyle \left (1 - \constbar{}^2\right)$ & $\displaystyle \left (1 - \frac{1}{m} \right )$ \\
	\hline
\end{tabular}
\end{table}
We use U$_\infty$ for the slower uniform rate to differentiate from U (tight uniform) for rate \eqref{eqapp:uniformtight}, and we use MR$_\infty$ for rate \eqref{maxResRate} to differentiate it from MR (tight)  rate \eqref{tightmaxres}.
 
For U$_\infty$, the rate follows straight from $\| A \|_{\infty,2} = \max_i \| a_i \| = \max_i \lambda_i = \lambda_{\max}$. For U, we note that the weighted matrix $\bar{A}:=D^{-1}A$ is simply the identity matrix. 
The NU rate uses that $\norm{A}_F^2 = \sum_i \lambda_i^2$.
For both MR$_\infty$ and MR, we have
\[
	\const{}^2 := \inf_{y \not = z } \frac{\| A (y - z ) \|_\infty^2}{\| y - z \|^2} = \inf_{\|w\| = 1} \| A w \|_\infty^2.
\]
Consider the equivalent problem
\begin{equation*}
\begin{split}
	\min_{w \in \R^m, \;y \in R} & \quad y \\
	\text{s.t.}   & \quad -y \le \lambda_i^2 w_i^2 \le y \text{ for all }i, \\
			& \quad \| w \| = 1,
\end{split}
\end{equation*}
From the first inequality, we get
\[
	- \frac{y}{\lambda_i^2} \le w_i^2 \le \frac{y}{\lambda_i^2} \quad \forall i \quad \Rightarrow \quad (w_i)^2 \le \frac{y}{\lambda_i^2} \quad \forall i.
\]
It follows that
\[
\norm{w}^2 = \sum_{i=1}^m w_i^2 \leq \sum_{i=1}^m \frac{y}{\lambda_i^2},
\]
which is equivalent to 
\[
y \geq \frac{\norm{w}^2}{\sum_{i=1}^m \frac{1}{\lambda_i^2}}.
\]
Because we are minimizing $y$ this must hold with equality at a solution, and because of the constraints $\norm{w} = 1$ we have
\[
	\const{}^2 = \left ( \sum_i \frac{1}{\lambda_i^2} \right )^{-1}.
\]
For the MR$_\infty$ rate, we divide $\const{}^2$ by the maximum eigenvalue squared. For the MR rate, we divide by the specific $\lambda_{i_k}^2$ corresponding to the row $i_k$ selected at iteration $k$.

For the MD rule, following the argument we did to derive $\const{}^2$ and using that $\bar{A} = {I}$ gives us
\[
	\constbar{}^2 = \frac{1}{m}.
\]

%


\section{Multiplicative Error}\label{app:multiplicative}

Suppose we have approximated the MR selection rule such that there is a multiplicative error in our selection of $i_k$, 
\[
	|a_{i_k}^T x^k - b_{i_k} | \ge \max_i |a_i^Tx^k - b_i| (1 - \epsilon_k),
\]
for some $\epsilon_k \in [0,1)$. In this scenario, we have 
\begin{align*}
    \| x^{k+1} - x^* \|^2
        &= \| x^k - x^* \|^2 - \frac{1}{\| a_{i_k} \|^2} \left(  \left | a_{i_k}^T x^k - b_{i_k} \right|^2\right) \\
        &\le \| x^k - x^* \|^2 - \frac{1}{\| a_{i_k} \|^2} \left( \max_i \left |a_i^T x^k - b_i \right| (1 - \epsilon_k )\right)^2 \\
	&= \| x^k - x^* \|^2 - \frac{(1 - \epsilon_k )^2}{\| a_{i_k} \|^2} \| A(x^k - x^*) \|^2_\infty \\
	&\le \bigg ( 1 - \frac{(1 - \epsilon_k )^2 \const{}^2}{\| a_{i_k} \|^2} \bigg ) \| x^k - x^* \|^2.
\end{align*}
We define a multiplicative approximation to the MD rule as an $i_k$ satisfying
\[
	\left |\frac{a_{i_k}^T x^k - b_{i_k}}{\|a_{i_k}\|} \right | \ge \max_i \left |\frac{a_i^Tx^k - b_i}{\|a_i\|} \right| (1 - \bar{\epsilon}_k),
\]
for some $\bar{\epsilon}_k \in [0,1)$. With such a rule we have
\begin{align*}
    \| x^{k+1} - x^* \|^2
        &= \| x^k - x^* \|^2 - \left(  \left | \frac{a_{i_k}^T x^k - b_{i_k}}{\|a_{i_k}\|} \right|^2\right) \\
        &\le \| x^k - x^* \|^2 - \left(  \max_i \left |\frac{a_i^Tx^k - b_i}{\|a_i\|} \right| (1 - \bar{\epsilon}_k)\right)^2 \\
        &= \| x^k - x^* \|^2 - (1 - \bar{\epsilon}_k )^2  \max_i ~\left |\frac{a_i^T(x^k - x^*)}{\|a_i\|} \right |^2 \\
	&= \| x^k - x^* \|^2 - (1 - \bar{\epsilon}_k )^2 \| D^{-1}A(x^k - x^*) \|^2_\infty \\
	&\le \bigg ( 1 - (1 - \bar{\epsilon}_k )^2 \constbar{}^2 \bigg ) \| x^k - x^* \|^2.
\end{align*}

\section{Additive Error}\label{app:additive}

Suppose we select $i_k$ using an approximate MR rule where
\[
	|a_{i_k}^T x^k - b_{i_k}|^2 \ge \max_i |a_i^Tx^k - b_i|^2 - \epsilon_k,
\]
for some $\epsilon_k \ge 0$.
Then we have the following convergence rate,
\begin{align*}
    \| x^{k+1} - x^* \|^2
        &= \| x^k - x^* \|^2 - \frac{1}{\| a_{i_k} \|^2} \left | a_{i_k}^T x^k - b_{i_k} \right |^2 \\
        &\le \| x^k - x^* \|^2 - \frac{1}{\| a_{i_k} \|^2}  \left (\max_i ~\left |a_i^T x^k - b_i \right|^2 - \epsilon_k  \right )\\
	&= \| x^k - x^* \|^2 - \frac{\| A(x^k - x^*) \|^2_\infty}{\| a_{i_k} \|^2}  + \frac{\epsilon_k}{\|a_{i_k}\|^2} \\
	&\le \bigg ( 1 - \frac{\const{}^2}{\| a_{i_k} \|^2} \bigg ) \| x^k - x^* \|^2 + \frac{\epsilon_k}{\|a_{i_k}\|^2}.
\end{align*}

For the MD rule with additive error, $i_k$ is selected such that
\[
	\left |\frac{a_{i_k}^T x^k - b_{i_k}}{\|a_{i_k} \|} \right|^2 \ge \max_i \left |\frac{a_i^T x^k - b_i}{\|a_i \|} \right|^2 - \bar{\epsilon}_k,
\]
for some $\bar{\epsilon_k} \ge 0$.
Then we have
\begin{align*}
	\| x^{k+1} - x^* \|^2
	&= \| x^k - x^* \|^2 -  \left |\frac{a_{i_k}^T x^k - b_{i_k}}{\|a_{i_k} \|} \right|^2\\
	&\le \| x^k - x^* \|^2 - \left( \max_i  \left | \frac{a_i^T x^k - b_i}{\|a_{i}\|} \right |^2 - \bar{\epsilon_k} \right ) \\
	&= \| x^k - x^* \|^2 -   \| D^{-1}A (x^k - x^*)  \|^2_\infty  + \bar{\epsilon_k} \\
	&\le \big ( 1 -   \constbar{}^2 \big ) \| x^k - x^* \|^2 + \bar{\epsilon_k}.
\end{align*}

\section{Comparison of Rates for the Maximum Distance Rule and the Randomized Kaczmarz via Johnson-Lindenstrauss Method}\label{app:73}

In \citet{eldar2011johnsonlindenstrauss}, the authors assume that the rows of $A$ are normalized and that we are dealing with a homogeneous system ($Ax = 0$), which is not particularly interesting since we can solve it in $O(1)$ by setting $x=0$. 
Their main convergence result is stated in Theorem 1. Note that {\em RKJL} stands for {\em Randomized Kaczmarz via Johnson-Lindenstrauss}, which is a hybrid technique using both random selection and an approximate MD rule using the dimensionality reduction technique of \citet{johnson1984lemma}. In their work they give the result below.

{\bf Theorem 1} ~{\em
	Fix an estimation $x^k$ and denote by $x^{k+1}$ and $x^{k+1}_{RK}$ the next estimations using the RKJL and the standard RK method, respectively. Define  $\gamma_j = |\langle a_j, x^k \rangle |^2$ and ordering these so that $\gamma_1 \ge \gamma_2 \ge \dots \ge \gamma_m$. Then, with $\delta$ being a constant affecting the error due to the JL approximation we have
\[
	\mathbb{E} \| x^{k+1} - x^* \|^2 \le \min \left [ \mathbb{E} \| x^{k+1}_{RK} - x \|^2 - \sum_{j = 1}^m \left ( p_j - \frac{1}{m} \right ) \gamma_j + 2 \delta, ~~\mathbb{E} \| x^{k+1}_{RK} - x^* \|^2 \right ],
\]
where
\[
	p_j = \begin{cases}
	\frac{\binom{m - j}{n - 1}}{\binom{m}{n}}, & j \le m - n + 1 \\
	0, & j > m - n + 1
	\end{cases}
\]
are non-negative values satisfying $\sum_{j = 1}^m p_j = 1$ and $p_1 \ge p_2 \ge \dots \ge p_m = 0$.
}

First, we simplify this bound. Applying the nonuniform random rate of \citet{strohmer2009} to the result of Theorem 1, we get
\begin{align*}
	&\mathbb{E} \left [\| x^{k+1} - x \| ^2 \right ] \\
	&\le \min \left [ \mathbb{E} \left [\| x^{k+1}_{RK} - x^* \|^2 \right ] - \sum_{j = 1}^m \left ( p_j - \frac{1}{m} \right ) \gamma_j + 2 \delta, ~~\mathbb{E} \left [ \| x^{k+1}_{RK} - x^* \|^2 \right ] \right ] \\
	&= \min \left [ \|x^k - x^* \|^2 - \frac{1}{\|A\|_F^2} \sum_{j = 1}^m \gamma_j - \sum_{j = 1}^m  p_j \gamma_j + \sum_{j = 1}^m \frac{1}{m} \gamma_j + 2 \delta, \|x^k - x^* \|^2 - \frac{1}{\|A\|_F^2} \sum_{j = 1}^m \gamma_j \right ] \\
	&= \min \left [ \|x^k - x^* \|^2   - \sum_{j = 1}^m  p_j \gamma_j + 2 \delta, ~~\|x^k - x^* \|^2 - \frac{1}{m} \sum_{j = 1}^m \gamma_j \right ], \numberthis \label{JLrate}
\end{align*}
where in the last line we use $\|A\|_F^2 = m$ for a matrix $A$ with normalized rows (in this case of normalized rows non-uniform selection  is simply uniform random selection). To compare this to our rate in the setting of an additive error, suppose we define $\epsilon_k$ such that the $i_k$ selected satisfies 
\[
	\gamma_{i_k} \ge \max_i \gamma_i - \bar{\epsilon}_k.
\]
Then, noting that $\|a_i\| = 1$ for all $i$, our convergence rate with additive error is based on the bound
\begin{align*}
    \| x^{k+1} - x^* \|^2
        &= \| x^k - x^* \|^2 - \gamma_{i_k}\\
        &\le \| x^k - x^* \|^2 - \max_i \gamma_i + \bar{\epsilon}_k. \numberthis \label{ourrateJL}
\end{align*}
Comparing the bounds \eqref{JLrate} and \eqref{ourrateJL}, we see that our MD bound
is always faster in the case of exact optimization ($\bar{\epsilon}_k = \delta = 0$), as the average and the weighted sum of the absolute inner products squared is less than the maximum inner product squared, $\max \{ \frac{1}{m} \sum_{j = 1}^m \gamma_j, \sum_{j = 1}^m  p_j \gamma_j \} \le \max_i \gamma_i$. If there is error present, then our rate is faster when
\[
	\max_i \gamma_i - \epsilon_k \ge \max \left \{ \frac{1}{m} \sum_{j = 1}^m \gamma_j, \sum_{j = 1}^m  p_j \gamma_j - 2 \delta \right \}.
\]
We note that even if our approximation is worse than the error resulting from the RKJL method, $\epsilon_k \ge 2 \delta$, it is possible that $\max_i \gamma_i$ is significantly larger than $\frac{1}{m} \sum_{j = 1}^m \gamma_j$ and $ \sum_{j = 1}^m  p_j \gamma_j$ and in this case our rate would be tighter. Further, our rate is more general as it does not specifically assume the Johnson-Lindenstrauss dimensionality reduction technique, that the rows of $A$ are normalized, or that the linear system is homogeneous. 

\section{Systems of Linear Inequalities}\label{app:inequalities}

Consider the system of linear equalities and inequalities,
\begin{equation}\label{eqapp:ineqsystem}
	\begin{cases}
		a_i^Tx \le b_i & (i \in I_\le) \\
		a_i^Tx = b_i & (i \in I_=).
	\end{cases}
\end{equation}
where the disjoint index sets $I_\le$ and $I_=$ partition the set $\{ 1,2,\dots,m\}$. As presented by \citet{leventhal2010constraints}, a generalization of the Kaczmarz algorithm that accommodates linear inequalities is given by
\begin{align*}
	\beta^k_{i_k} & =
	\begin{cases}
		(a_{i_k}^T x^k - b_{i_k})^+ & ({i_k} \in I_\le) \\
		~a_{i_k}^T x^k - b_{i_k} & ({i_k} \in I_=),
	\end{cases} \\
	x^{k+1} &= x^k - \frac{\beta^k_{i_k}}{\| a_{i_k} \|^2} a_{i_k},
\end{align*}
where for $x \in \R^n$ we define $x^+$  element-wise by
\[
	(x^+)_i = \max \{ x_i,0 \}.
\]

This leads to the following generalization of the MR and MD rules, respectively,
\begin{equation}\label{genMax}
	i_k = \max \left | \beta^k_i \right | = \| \beta^k \|_\infty , \quad \text{and} \quad i_k = \max \left | \frac{\beta^k_i}{\|a_i \|} \right | = \| D^{-1} \beta^k \|_\infty.
\end{equation}
Unlike for equalities where the Kaczmarz method converges to the projection of the initial iterate $x^0$ onto the intersection of the constraints, for inequalities we can only guarantee that the Kaczmarz method converges to a point in the feasible set. Thus, in convergence rates involving inequalities it is standard to use a bound for the distance from the current iterate $x^k$ to the feasible region,
\[
	d(x,S) = \min_{z \in S} \| x - z \|_2 = \| x - P_S(x) \|_2,
\]
where $P_S(x)$ is the projection of $x$ onto the feasible set $S$. 

Following closely the arguments of \citet{leventhal2010constraints} for systems of inequalities, we next give the following result which they credit to~ \citet{hoffman1952}.
%

\begin{theorem}\label{thm:hoffman}
Let \eqref{eqapp:ineqsystem} be a consistent system of linear equalities and inequalities,
then there exists a constant $\const{}$ such that
\[
	x \in \R^n \text{ and } S \not = \emptyset \quad \Rightarrow \quad d(x,S) \le \frac{1}{\const{}} \| e(Ax-b) \|_\infty,
\]
where $S$ is the set of feasible solutions and
where the function $e:\R^m \mapsto \R^m$ is defined by 
\[
	e(y)_i = 
	\begin{cases}
		y_i^+ & (i \in I_\le) \\
		y_i & (i \in I_=).
	\end{cases}
\]
\end{theorem}
From \citet{leventhal2010constraints}, combining both cases ($i_k \in I_\le$ or $i_k \in I_=$), the following relationship holds with respect to the distance measure $d(x,S)$,
\begin{equation}\label{ineqproof}
	d(x^{k+1}, S)^2 \le d(x^k, S) - \frac{e(Ax^k - b)_{i_k}^2}{\| a_{i_k} \|^2}.
\end{equation}
Following from this bound and Theorem \ref{thm:hoffman}, it is straightforward to derive analogous results for all greedy selection rates derived in the paper. For example, if we select $i_k$ according to the generalized MR rule \eqref{genMax} then the analogous tight rate for the MR rule is given by
\begin{align*}
	d(x^{k+1}, S)^2 
	&\le d(x^k, S)^2 - \frac{e(Ax^k - b)_{i_k}^2}{\| a_{i_k} \|^2} \\
	&= d(x^k, S)^2 - \frac{ \| \beta^k \|^2_\infty}{\| a_{i_k} \|^2} \\
	&\le \left ( 1- \frac{\const{}^2}{\| a_{i_k} \|^2} \right ) d(x^k, S)^2.
\end{align*}

\section{Multi-Step Maximum Residual Bound}\label{app:multi-MR}

Recall the MR rate~\eqref{rate:tightmaxres},
\[
	\| x^k - x^* \|^2 \le \prod_{i = 1}^k\bigg ( 1 - \frac{\const{}^2}{\| a_{i_{k}} \|^2} \bigg ) \| x^0 - x^* \|^2.
\]
In the worst case this is no faster than the MR$_\infty$ rate since we may have $\norm{a_{i_k}} = \norm{A}_{\infty,2}$ for all $i$. However, this rate is faster if we have $\norm{a_{i_k}} < \norm{A}_{\infty,2}$ for any $i$. In this section we derive a  bound that will typically be much tighter than MR$_\infty$ by considering the sequence of $\norm{a_{i_k}}$ values that are possible for problems with a sparse orthogonality graph.
To derive an upper bound, we solve the problem below which was first introduced in~\citet{nutini2015}.
\begin{problem}{{\bf1}.}\label{prob:1}
{\em 
We are given a graph $G = (V,E)$ with $E \neq \emptyset$, a weight $M_i$ associated with each node $i$, and an iteration number $k$. Choose a sequence $\{ i_t \}_{t = 1}^k$ that maximizes the sum of the weights $M_{i_t}$ subject to the following constraint: after each time node $i$ has been chosen, it cannot be chosen again until after a neighbour of node $i$ has been chosen.
}	
\end{problem}

 To map this problem to the problem of showing that the $\norm{a_{i_k}}$ values are small when we use the MR rule, we the weights $M_{i_k} =  \log \left(1 - \frac{\const{}^2}{\| a_{i_{k}} \|^2} \right )$. The constraint in Problem 1 arises because the MR rule cannot choose $i_k$ on any future iteration until after a neighbour of it is selected in the orthogonality graph.  \citet{nutini2015} give a bound on the solution of this problem for the case of chain-structured graphs, but we have now derived the result for the case of general graphs. In particular, it can be shown that an asymptotically optimal sequence of weights is given by repeatedly visiting the star subgraph of the original graph with the maximum average weight. Using this result and the mapping to the Kaczmarz problem yields the bound stated in the main paper. The proof of the asymptotic optimality of the maximum-average-weight star-structured subgraph is highly non-trivial, and is contained in the MSc thesis of the second author~\citep[][Theorem~2.1]{behrooz2016}
  
\section{Faster Randomized Kaczmarz Using the Orthogonality Graph of $A$}\label{app:fasterrandom}

In order for the adaptive methods to be efficient, we must be able to efficiently update the set of selectable nodes at each iteration. 
To do this we use a tree structure that keeps track of the number of selectable children in the tree (for uniform random selection) or the cumulative sums of the selectable row norms of $A$ (for non-uniform random selection). A similar structure is used in the non-uniform sampling code of~\citet{schmidt2013sag}.

Recall that the standard inverse-transform approach approach to sampling from a non-uniform discrete probability distribution over $m$ variables:
\begin{enumerate}
	\item Compute the cumulative probabilities, $c_i = \sum_{j=1}^i p_j$ for each $i$ from $1$ to $m$.
	\item Generate a random number $u$ uniformly distributed over $[0,1]$.
	\item Return the smallest $i$ such that $c_i \geq u$.
\end{enumerate}
We can compute all $m$ values of  $c_i$ in Step 1 at a cost of $O(m)$ by maintaining the running sum. We'll assume that Step 2 costs $O(1)$ and we can implement Step 3 in $O(\log(m))$ using a binary search. If we are sampling from a fixed distribution, then we only need to perform Step 1 once and from that point we can generate samples from the distribution at a cost of $O(\log(m))$.

In the adaptive randomized selection rules, the probabilities $p_j$ change at each iteration and hence the $c_i$ values also change. This means we can't skip Step 1 as we can for fixed probabilities. However, if the orthogonality graph is sparse then it's still possible to efficiently implement these strategies. To do this, we consider a binary tree-structure that has the probabilities $p_j$ as leaf nodes while each internal node is the {\em sum} of its two descendants (and thus the root node has a value of 1). Given this structure, we can find the smallest $c_i \geq u$ in $O(\log(m))$ by traversing the tree. Further, if we update one of the $p_j$ values then we can update this data structure in $O(\log(m))$ time since this only requires changing one node at each depth of the tree. If each node has at most $g$ neighbours in the orthogonality graph, then we need to update $g$ probabilities in the binary tree, leading to a cost of $O(g \log(m))$ to update the tree structure on each iteration. 
 
Note that the above structure can be modified to work with \emph{unnormalized probabilities} at the leaf nodes, since the root node will contain the normalizing constant required to make these unnormalized probabilities into a valid probability mass function. Using this, we can implement the adaptive uniform method by setting the leaf nodes to $1$ for selectable nodes and $0$ for non-selectable nodes. To implement the adaptive non-uniform method, we set the leaf nodes to $0$ for non-selectable nodes and $\norm{a_i}^2$ for selectable nodes.

\section{Additional Experiments}\label{app:numerics}

\subsection*{Formulating the Semi-Supervised Label Propagation Problem as a Linear System}

Our third experiment solves a label propagation problem for semi-supervised learning in the `two moons' dataset \citep{zhou2004learning}. 
We use a variant of the quadratic labelling criterion of \citet{bengio2006},
\[
	\min_{y_i \in S' } \quad f(y) \equiv \frac{1}{2} \sum_{i=1}^n \sum_{j=1}^n w_{ij} (y_i - y_j)^2,
\]
where $y$ is our label vector (each $y_i$ can take one of 2 values), $S$ is the set of labels that we do know, $S'$ is the set of labels that we do not know and $w_{ij} \ge 0$ are the weights assigned to each $y_i$ describing how strongly we want the labels $y_i$ and $y_j$ to be similar. We assume without loss of generality that $w_{ii} = 0$ (since it doesn't affect the objective) and that $w_{ij} = w_{ji}$ for all $i,j$ because by the symmetry in the objective the model only depends on these terms through $(w_{ij} + w_{ji})$. We can express this quadratic problem as a linear system that is consistent by construction. In other words, we can define $A$ and $b$ such that
\[
	\nabla f(y) = 0 \iff Ay = b, \quad \text{with } y \in S'.
\]
Differentiating $f$ with respect to some $y_k \in S'$, we have
\begin{align*}
	\nabla_k f(y) 
	&= \underbrace{\sum_{j \not = k} w_{kj} (y_k - y_j)}_{i = k, ~j \not = k} 
	- \underbrace{\sum_{i \not = k} w_{ik} (y_i - y_k)}_{i \not = k, ~j  = k}
	+ \underbrace{\sum_{i = k} w_{kk} (y_k - y_k)}_{i  = k, ~j  = k} \\
	&= \sum_{i = 1}^n w_{ki} (y_k - y_i)
	- \sum_{i = 1}^n w_{ik} (y_i - y_k) \\
	&= 2\sum_{i =1}^n w_{ki} y_k - 2\sum_{i = 1}^n w_{ki} y_i.  
\end{align*}
Setting this equal to zero and splitting the summation over $S$ and $S'$ separately, we have
\begin{align*}
	 \sum_{i = 1}^n w_{ki} y_k  - \sum_{i \in S'} w_{ki}y_i
	&= \sum_{i \in S} w_{ki} y_i.
\end{align*}
Assuming the elements of $S'$ form the first $|S'|$ elements of the matrix $A$, the above formulation yields the  $|S'| \times |S'|$ matrix with entries
\[
	A_{k,i} = 
	\begin{cases}
		\sum_{j=1}^n w_{kj} & \text{if } i = k, \\
		-w_{ki} & \text{if } i \not = k,
	\end{cases}	
\]
where $k$ and $i \in S'$ and
\[
	b_k  = \sum_{i \in S} w_{ki}y_i.
\]

\subsection*{Time vs. Squared Error and Distance}

\begin{figure*}
\centering
\includegraphics[width=.49\textwidth]{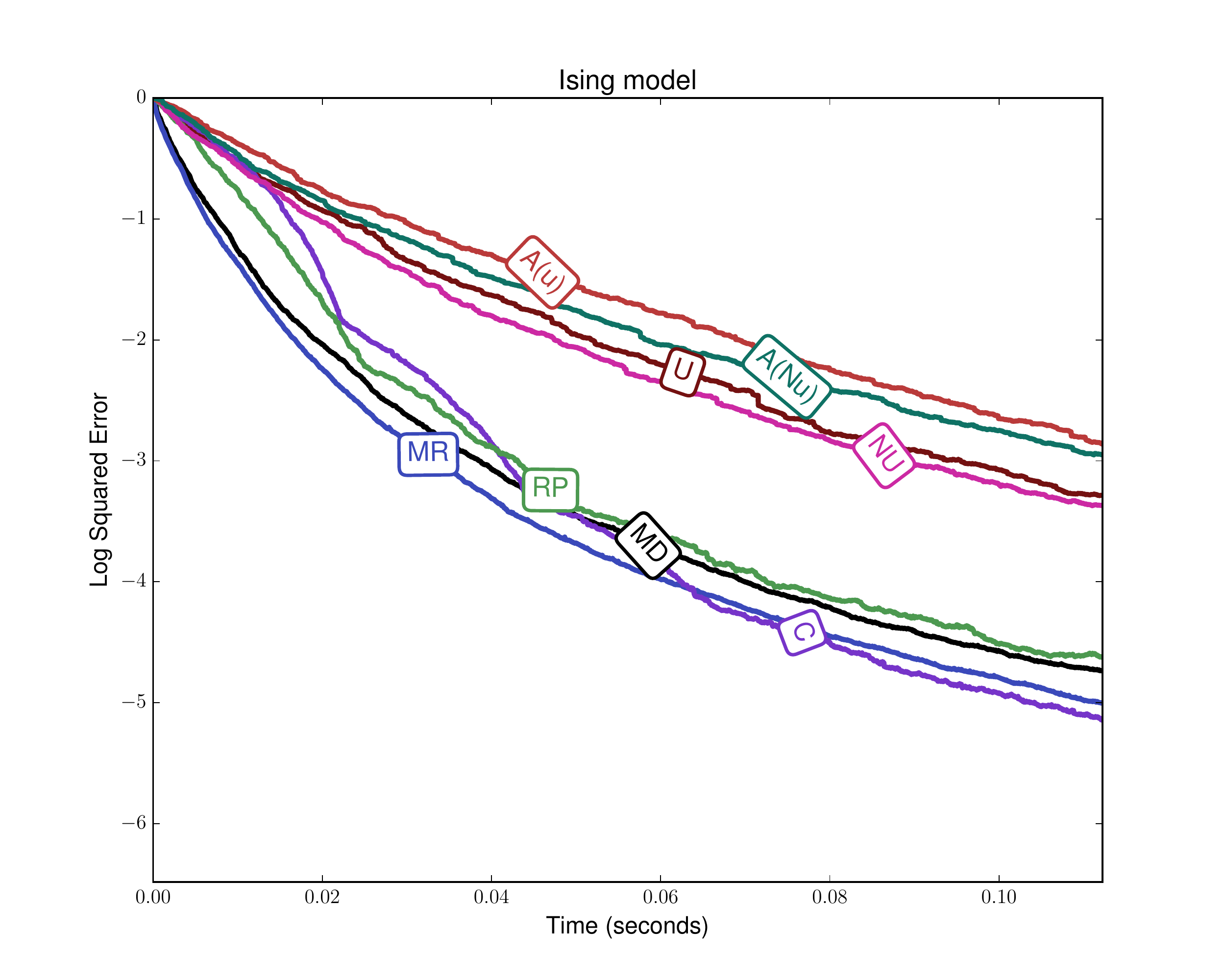}
\includegraphics[width=.49\textwidth]{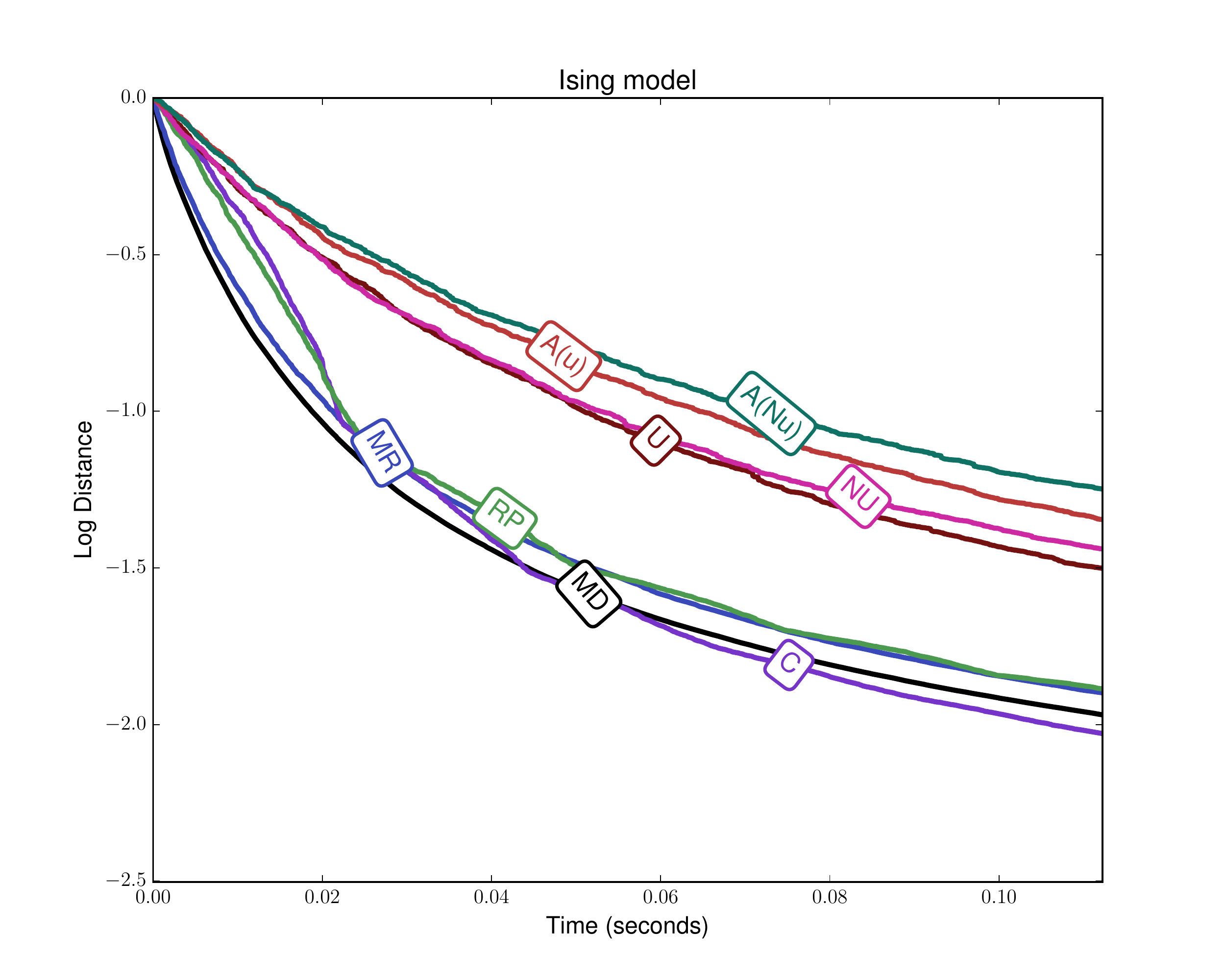}\\ 
\includegraphics[width=.49\textwidth]{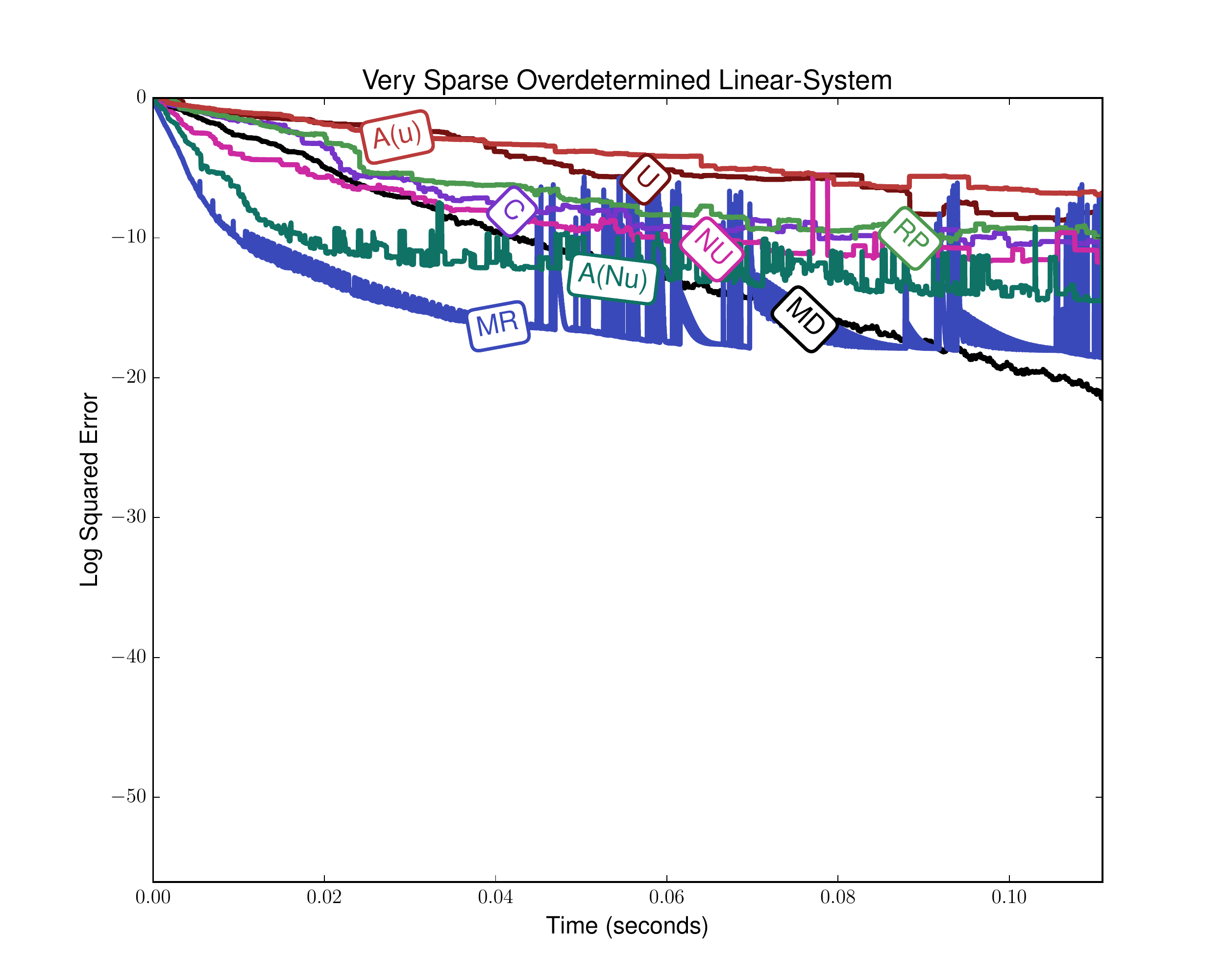}
\includegraphics[width=.49\textwidth]{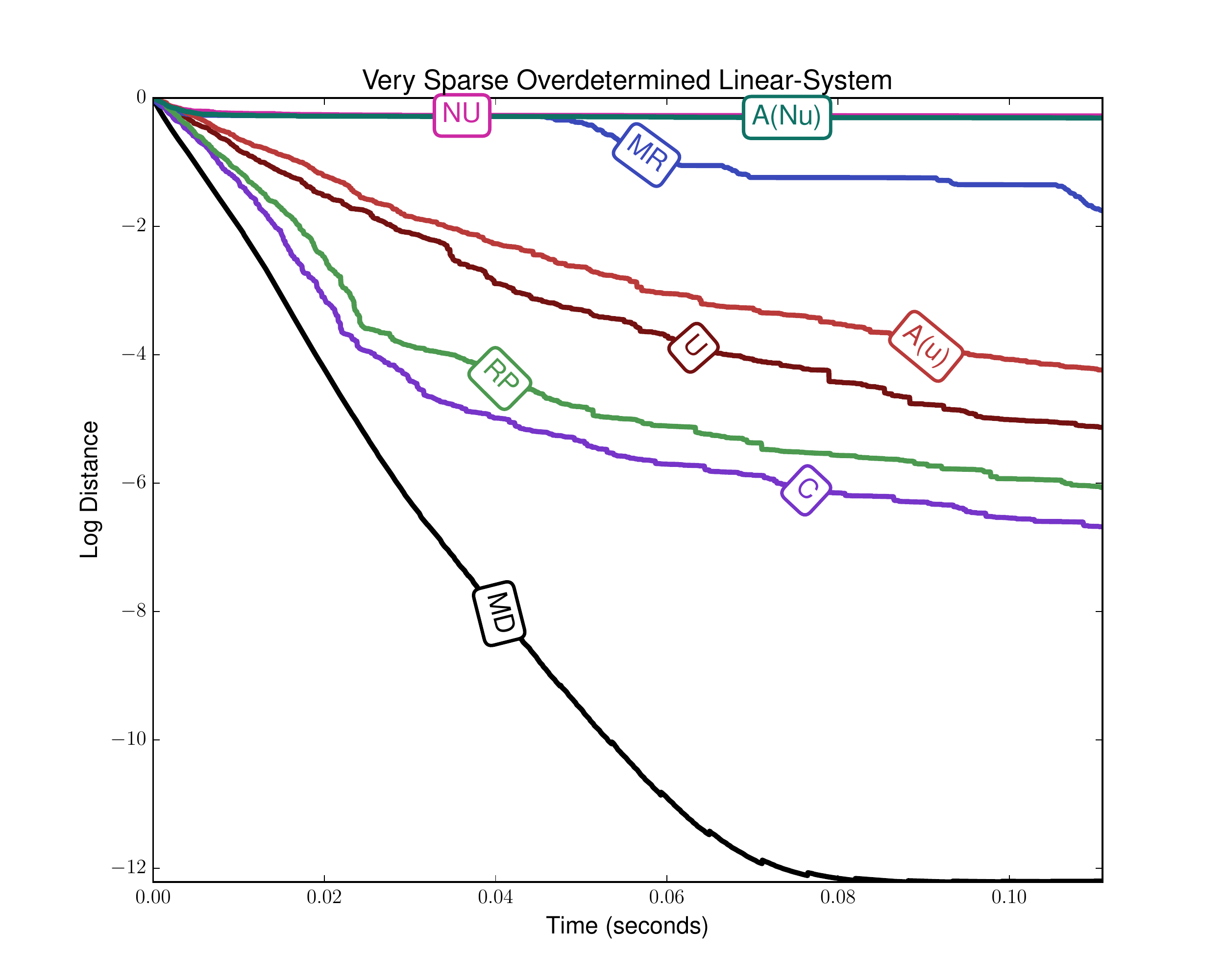}\\ 
\includegraphics[width=.49\textwidth]{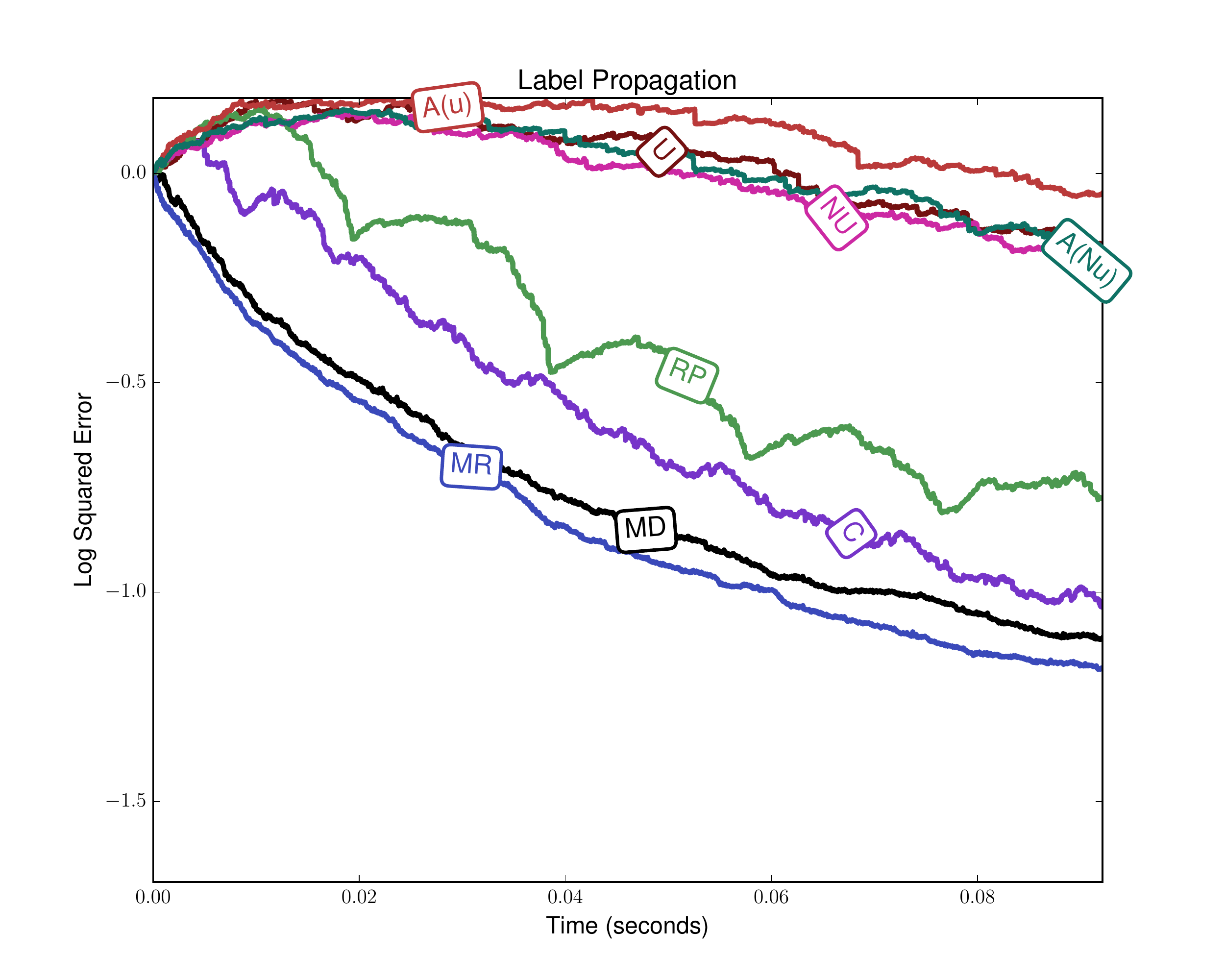}  
\includegraphics[width=.49\textwidth]{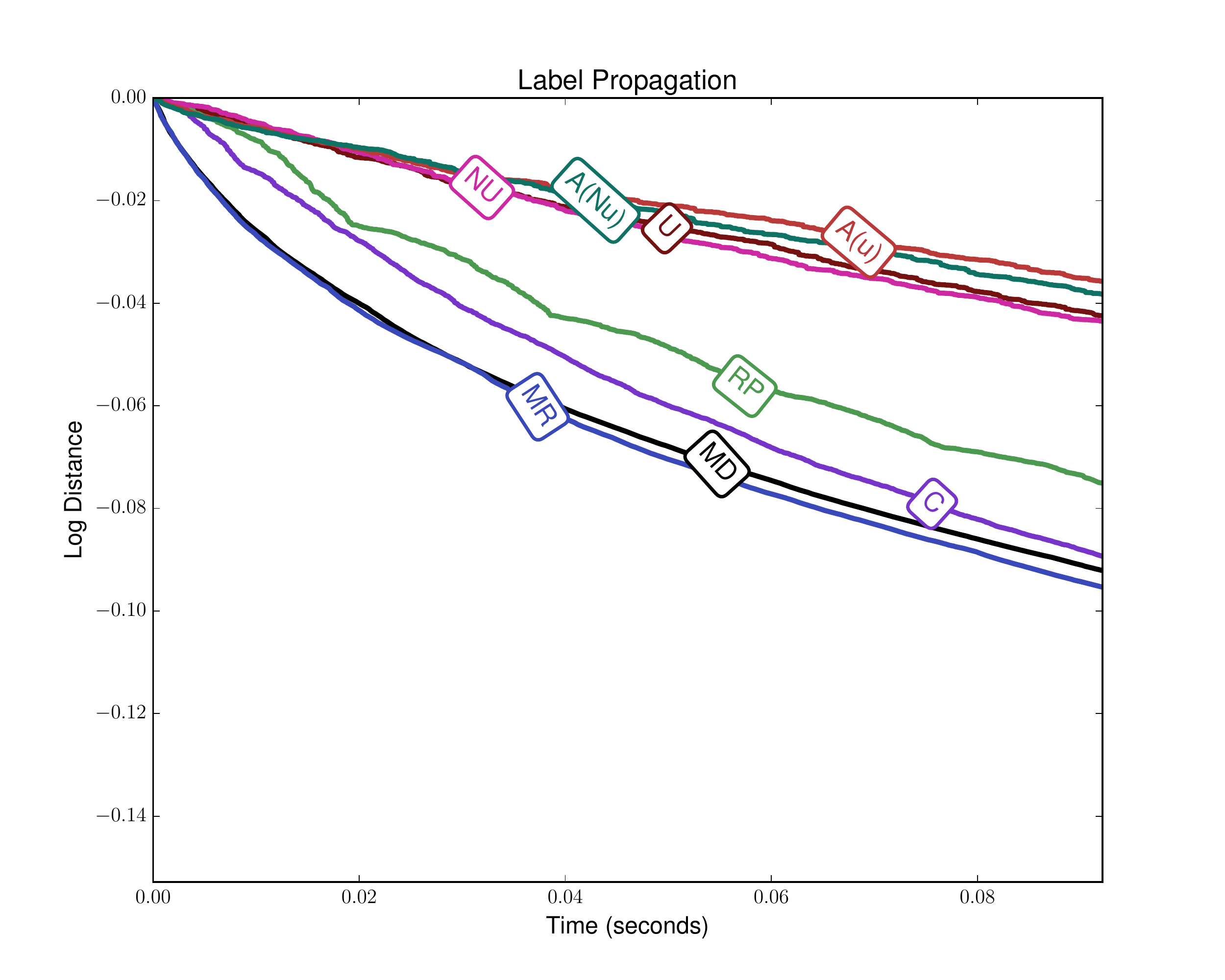} 
\caption{Runtime Comparisons of Kaczmarz Selection Rule.}
\label{fig:exp_time}
\end{figure*}

Figure \ref{fig:exp_time} compares the runtime results for our 3 experiments from the main paper using both squared error and distance (we made a reasonable effort to make the implementations of all methods as efficient as possible). We see that in the first experiment the greedy selection rules do not translate into gains in terms of runtime over the cyclic methods due to their higher iteration cost (though they still outperform random methods), while in the second and third experiments the greedy rules are still superior in terms of runtime.


\subsection*{Hybrid Methods}

\begin{figure*}[!ht]
\centering
\includegraphics[width=.49\textwidth]{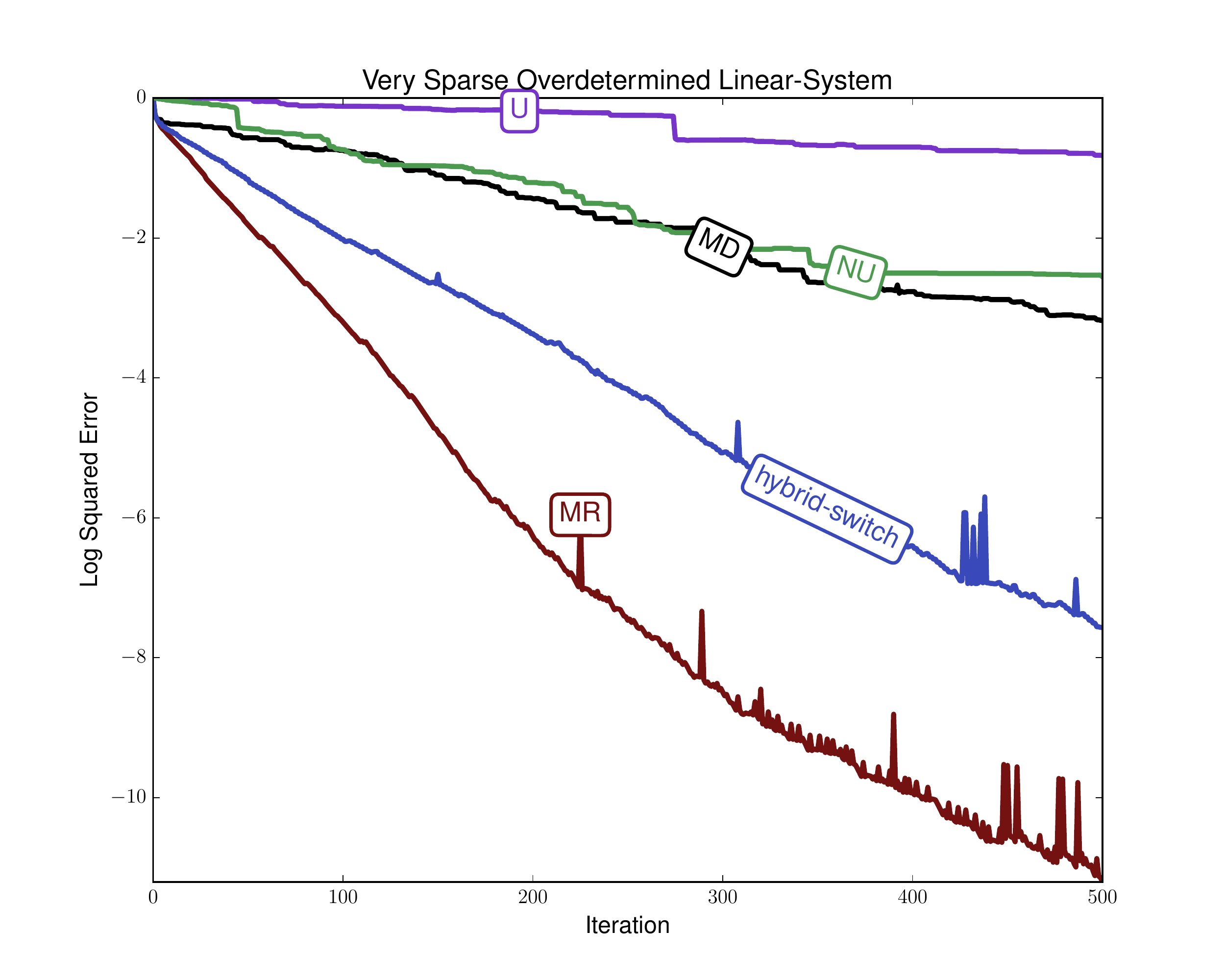}
\includegraphics[width=.49\textwidth]{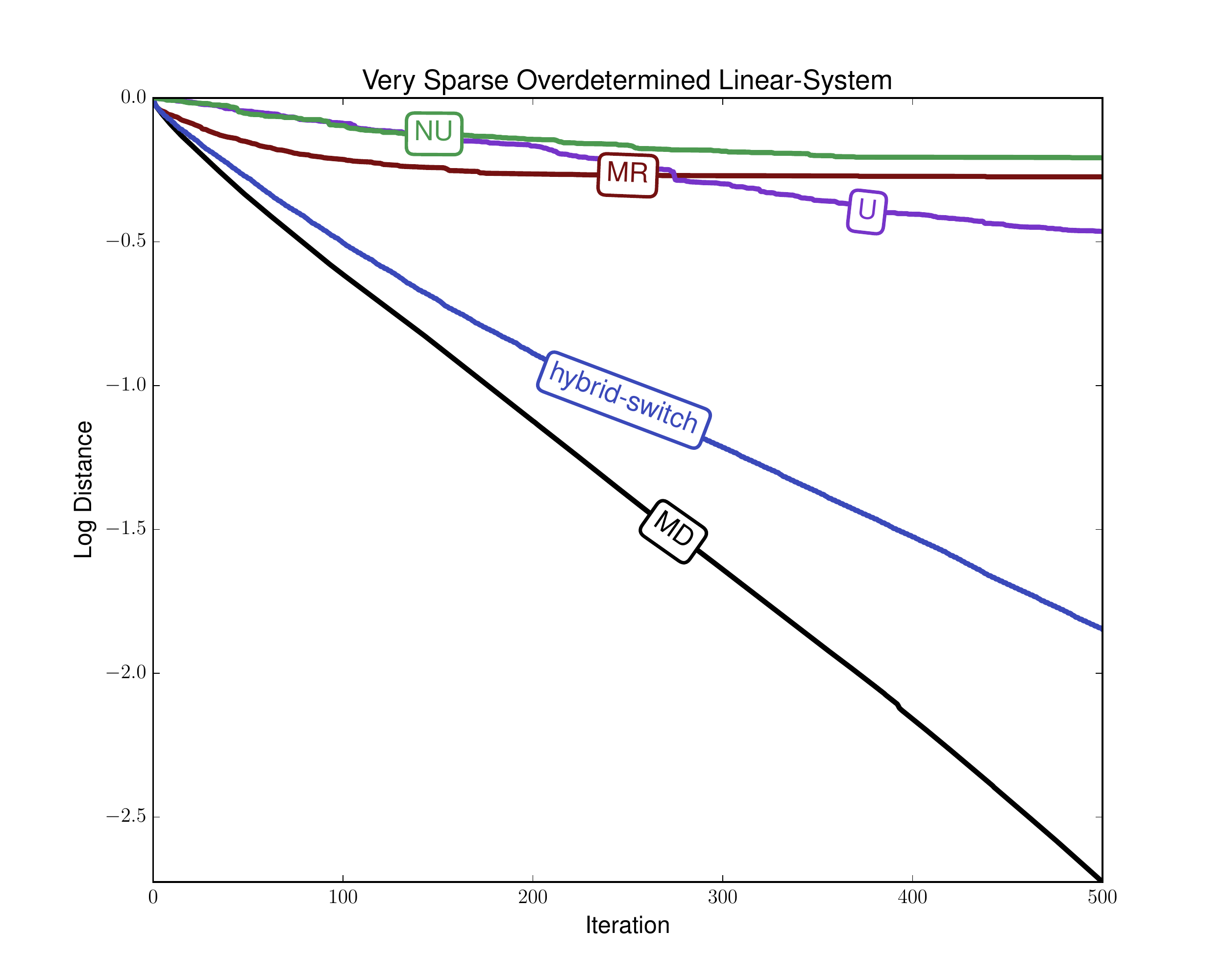}
\caption{Comparison of MR, MD and Hybrid Method for Very Sparse Dataset.}
\label{fig:hybrid}
\end{figure*}

For the very sparse overdetermined dataset, we see very different performances between the MR and MD rules with respect to squared error and distance. We see that the MR rule outperforms the MD rule in the beginning with respect to squared-error and the MD rule outperforms the MR rule significantly with respect to distance. These observations align with the respective definitions of each greedy rule. However, if we want a method that converges well with respect to {\em both} of these objectives, then we could consider `hybrid' greedy rule. For example, we could simply alternate between using the MR rule and the MD rule. As we see in Figure \ref{fig:hybrid}, this approach simultaneously exploits the convergence of the MR rule in terms of squared error and the MD rule in terms of distance to the solution.
However, computationally this approach requires the maintenance of two max-heap structures.

\subsubsection*{References}

{
\bibliography{bib}

\begin{thebibliography}{39}
\providecommand{\natexlab}[1]{#1}
\providecommand{\url}[1]{\texttt{#1}}
\expandafter\ifx\csname urlstyle\endcsname\relax
  \providecommand{\doi}[1]{doi: #1}\else
  \providecommand{\doi}{doi: \begingroup \urlstyle{rm}\Url}\fi

\bibitem[Bengio et~al.(2006)Bengio, Delalleau, and Le~Roux]{bengio2006}
Y.~Bengio, O.~Delalleau, and N.~Le~Roux.
\newblock Label propagation and quadratic criterion.
\newblock In O.~Chapelle, B.~Sch\"olkopf, and A.~Zien, editors,
  \emph{Semi-Supervised Learning}, chapter~11, pages 193--216. MIT Press, 2006.

\bibitem[Censor(1981)]{censor1981}
Y.~Censor.
\newblock Row-action methods for huge and sparse systems and their
  applications.
\newblock \emph{SIAM Rev.}, 23\penalty0 (4):\penalty0 444--466, 1981.

\bibitem[Censor et~al.(1983)Censor, Eggermont, and Gordon]{censor1983}
Y.~Censor, P.~B. Eggermont, and D.~Gordon.
\newblock Strong underrelaxation in {K}aczmarz's method for inconsistent
  systems.
\newblock \emph{Numer. Math.}, 41:\penalty0 83--92, 1983.

\bibitem[Censor et~al.(2009)Censor, Herman, and Jiang]{censor2009notesSV}
Y.~Censor, G.~T. Herman, and M.~Jiang.
\newblock A note on the behaviour of the randomized {K}aczmarz algorithm of
  {S}trohmer and {V}ershynin.
\newblock \emph{J. Fourier Anal. Appl.}, 15:\penalty0 431--436, 2009.

\bibitem[Deutsch(1985)]{deutsch1985}
F.~Deutsch.
\newblock Rate of convergence of the method of alternating projections.
\newblock \emph{Internat. Schriftenreihe Numer. Math.}, 72:\penalty0 96--107,
  1985.

\bibitem[Deutsch and Hundal(1997)]{deutsch1997}
F.~Deutsch and H.~Hundal.
\newblock The rate of convergence for the method of alternating projections,
  {II}.
\newblock \emph{J. Math. Anal. Appl.}, 205:\penalty0 381--405, 1997.

\bibitem[Eldar and Needell(2011)]{eldar2011johnsonlindenstrauss}
Y.~C. Eldar and D.~Needell.
\newblock Acceleration of randomized {K}aczmarz methods via the
  {J}ohnson-{L}indenstrauss {L}emma.
\newblock \emph{Numer. Algor.}, 58:\penalty0 163--177, 2011.

\bibitem[Feichtinger et~al.(1992)Feichtinger, Cenker, Mayer, Steier, and
  Strohmer]{cenker1992POCS}
H.~G. Feichtinger, C.~Cenker, M.~Mayer, H.~Steier, and T.~Strohmer.
\newblock New variants of the {POCS} method using affine subspaces of finite
  codimension with applications to irregular sampling.
\newblock \emph{SPIE: VCIP}, pages 299--310, 1992.

\bibitem[Gal\'antai(2005)]{galantai2005}
A.~Gal\'antai.
\newblock On the rate of convergence of the alternating projection method in
  finite dimensional spaces.
\newblock \emph{J. Math. Anal. Appl.}, 310:\penalty0 30--44, 2005.

\bibitem[Gordon et~al.(1970)Gordon, Bender, and Herman]{gordon1970ART}
R.~Gordon, R.~Bender, and G.~T. Herman.
\newblock Algebraic {R}econstruction {T}echniques ({ART}) for three-dimensional
  electron microscopy and x-ray photography.
\newblock \emph{J. Theor. Biol.}, 29\penalty0 (3):\penalty0 471--481, 1970.

\bibitem[Gower and Richt\'{a}rik(2015)]{gower2015}
R.~M. Gower and P.~Richt\'{a}rik.
\newblock Randomized iterative methods for linear systems.
\newblock \emph{SIAM J. Matrix Anal. Appl.}, 36\penalty0 (4):\penalty0
  1660--1690, 2015.

\bibitem[Griebel and Oswald(2012)]{griebel2012}
M.~Griebel and P.~Oswald.
\newblock Greedy and randomized versions of the multiplicative {S}chwartz
  method.
\newblock \emph{Lin. Alg. Appl.}, 437:\penalty0 1596--1610, 2012.

\bibitem[Hanke and Niethammer(1990)]{hanke1990}
M.~Hanke and W.~Niethammer.
\newblock On the acceleration of {K}aczmarz's method for inconsistent linear
  systems.
\newblock \emph{Lin. Alg. Appl.}, 130:\penalty0 83--98, 1990.

\bibitem[Herman and Meyer(1993)]{herman1993}
G.~T. Herman and L.~B. Meyer.
\newblock Algebraic reconstruction techniques can be made computationally
  efficient.
\newblock \emph{IEEE Trans. Medical Imaging}, 12\penalty0 (3):\penalty0
  600--609, 1993.

\bibitem[Hoffman(1952)]{hoffman1952}
A.~J. Hoffman.
\newblock On approximate solutions of systems of linear inequalities.
\newblock \emph{J. Res. Nat. Bur. Stand.}, 49\penalty0 (4):\penalty0 263--265,
  1952.

\bibitem[Johnson and Lindenstrauss(1984)]{johnson1984lemma}
W.~B. Johnson and J.~Lindenstrauss.
\newblock Extensions of {L}ipchitz mappings into a {H}ilbert space.
\newblock \emph{Contemp. Math.}, 26:\penalty0 189--206, 1984.

\bibitem[Kaczmarz(1937)]{kaczmarz1937}
S.~Kaczmarz.
\newblock Angen\"aherte {A}ufl\"osung von {S}ystemen linearer {G}leichungen,
  {B}ulletin {I}nternational de l'{A}cad\'emie {P}olonaise des {S}ciences et
  des {L}etters.
\newblock \emph{Classe des Sciences Math\'ematiques et Naturelles. S\'erie A,
  Sciences Math\'ematiques}, 35:\penalty0 355--357, 1937.

\bibitem[Karimi et~al.(2016)Karimi, Nutini, and Schmidt]{karimi2016}
H.~Karimi, J.~Nutini, and M.~Schmidt.
\newblock Linear convergence of gradient and proximal-gradient methods under
  the polyak-\l{}ojasiewicz condition.
\newblock \emph{European Conference on Machine Learning}, 2016.

\bibitem[Lee and Sidford(2013)]{lee2013}
Y.~T. Lee and A.~Sidford.
\newblock Efficient accelerated coordinate descent methods and faster
  algorithms for solving linear systems.
\newblock \emph{arXiv:1305.1922v1}, 2013.

\bibitem[Leventhal and Lewis(2010)]{leventhal2010constraints}
L.~Leventhal and A.~S. Lewis.
\newblock Randomized methods for linear constraints: convergence rates and
  conditioning.
\newblock \emph{Math. Oper. Res.}, 35\penalty0 (3):\penalty0 641--654, 2010.

\bibitem[Liu and Wright(2014)]{liu2014accelerated}
J.~Liu and S.~J. Wright.
\newblock An accelerated randomized {K}aczmarz method.
\newblock \emph{arXiv:1310.2887v2}, 2014.

\bibitem[Ma et~al.(2015)Ma, Needell, and Ramdas]{ma2015}
A.~Ma, D.~Needell, and A.~Ramdas.
\newblock Convergence properties of the randomized extended {G}auss-{S}eidel
  and {K}aczmarz methods.
\newblock \emph{arXiv:1503.08235v2}, 2015.

\bibitem[Needell(2010)]{needell2010noisy}
D.~Needell.
\newblock Randomized {K}aczmarz solver for noisy linear systems.
\newblock \emph{BIT Numer. Math.}, 50:\penalty0 395--403, 2010.

\bibitem[Needell and Tropp(2014)]{needell2012}
D.~Needell and J.~A. Tropp.
\newblock Paved with good intentions: {A}nalysis of a randomized block
  {K}aczmarz method.
\newblock \emph{Lin. Alg. Appl.}, 441:\penalty0 199--221, 2014.

\bibitem[Needell et~al.(2015)Needell, Srebro, and Ward]{needell2014sgd}
D.~Needell, N.~Srebro, and R.~Ward.
\newblock Stochastic gradient descent and the randomized {K}aczmarz algorithm.
\newblock \emph{arXiv:1310.5715v5}, 2015.

\bibitem[Novikoff(1962)]{novikoff1963convergence}
A.~B.~J. Novikoff.
\newblock On convergence proofs for perceptrons.
\newblock \emph{Symp. Math. Theory Automata}, 12:\penalty0 615--622, 1962.

\bibitem[Nutini et~al.(2015)Nutini, Schmidt, Laradji, Friedlander, and
  Koepke]{nutini2015}
J.~Nutini, M.~Schmidt, I.~H. Laradji, M.~Friedlander, and H.~Koepke.
\newblock Coordinate descent converges faster with the {G}auss-{S}outhwell rule
  than random selection.
\newblock \emph{ICML}, 2015.

\bibitem[Oswald and Zhou(2015)]{oswald2015}
P.~Oswald and W.~Zhou.
\newblock Convergence analysis for {K}aczmarz-type methods in a {H}ilbert space
  framework.
\newblock \emph{Lin. Alg. Appl.}, 478:\penalty0 131--161, 2015.

\bibitem[Rue and Held(2005)]{rue2005}
H.~Rue and L.~Held.
\newblock \emph{Gaussian Markov Random Fields: Theory and Applications}.
\newblock CRC Press, 2005.

\bibitem[Schmidt et~al.(2013)Schmidt, {Le Roux}, and Bach]{schmidt2013sag}
M.~Schmidt, N.~{Le Roux}, and F.~Bach.
\newblock Minimizing finite sums with the stochastic average gradient.
\newblock \emph{arXiv preprint}, 2013.

\bibitem[Sepehry(2016)]{behrooz2016}
B.~Sepehry.
\newblock Finding a maximum weight sequence with dependency constraints.
\newblock Master's thesis, University of British Columbia, 2016.

\bibitem[Strohmer and Vershynin(2009)]{strohmer2009}
T.~Strohmer and R.~Vershynin.
\newblock A randomized {K}aczmarz algorithm with exponential convergence.
\newblock \emph{J. Fourier Anal. Appl.}, 15:\penalty0 262--278, 2009.

\bibitem[Tanabe(1971)]{tanabe1971}
K.~Tanabe.
\newblock Projection method for solving a singular system of linear equations
  and its applications.
\newblock \emph{Numer. Math.}, 17:\penalty0 203--214, 1971.

\bibitem[van~der Sluis(1969)]{vandersluis1969}
A.~van~der Sluis.
\newblock Condition numbers and equilibrium of matrices.
\newblock \emph{Numer. Math.}, 14:\penalty0 14--23, 1969.

\bibitem[Vishnoi(2013)]{vishnoi2012laplacian}
N.~K. Vishnoi.
\newblock ${L}x = b$ {L}aplacian solvers and their algorithmic applications.
\newblock \emph{Found. Trends Theoretical Computer Science}, 8\penalty0
  (1-2):\penalty0 1--141, 2013.

\bibitem[Whitney and Meany(1967)]{whitney1967}
T.~Whitney and R.~Meany.
\newblock Two algorithms related to the method of steepest descent.
\newblock \emph{SIAM J. Numer. Anal.}, 4\penalty0 (1):\penalty0 109--118, 1967.

\bibitem[Wright(2015)]{wright2015}
S.~J. Wright.
\newblock Coordinate descent algorithms.
\newblock \emph{arXiv:1502.04759v1}, 2015.

\bibitem[Zhou et~al.(2004)Zhou, Bousquet, Lal, Weston, and
  Sch{\"o}lkopf]{zhou2004learning}
D.~Zhou, O.~Bousquet, T.~N. Lal, J.~Weston, and B.~Sch{\"o}lkopf.
\newblock Learning with local and global consistency.
\newblock \emph{NIPS}, 2004.

\bibitem[Zouzias and Freris(2013)]{zouzais2012extended}
A.~Zouzias and N.~M. Freris.
\newblock Randomized extended {K}aczmarz for solving least-squares.
\newblock \emph{arXiv:1205.5770v3}, 2013.

\end{thebibliography}
\bibliographystyle{abbrvnat}
}

\end{document}